\documentclass[12pt]{amsart}

\newtheorem{lem}{Lemma}[section]
\newtheorem{defi}[lem]{Definition}
\newtheorem{prop}[lem]{Proposition}
\newtheorem{rema}[lem]{Remark}
\newtheorem{theo}[lem]{Theorem}
\newtheorem{coro}[lem]{Corollary}
\newtheorem{exam}[lem]{Example}

\begin{document}
\parindent0em
\title{One parameter fixed point theory and gradient flows of closed 1-forms}
\author{D. Sch\"utz}
\address{Department of Mathematics, SUNY Binghamton, Binghamton, NY 13902-6000}
\email{dirk@math.binghamton.edu}
\begin{abstract}
We use the one-parameter fixed point theory of Geoghegan and Nicas to get information about the closed
orbit structure of transverse gradient flows of closed 1-forms on a closed manifold $M$.
We define a noncommutative zeta function in an object related to the first Hochschild homology
group of the Novikov ring associated to the 1-form and relate it to the torsion of a natural chain
homotopy equivalence between the Novikov complex and a completed simplicial complex of $\tilde{M}$,
the universal cover of $M$.
\end{abstract}
\maketitle
\section{Introduction}
Let $\omega$ be a closed 1-form on a closed connected smooth manifold $M$. There is a corresponding
homomorphism $\xi:G\to\mathbb{R}$, where $G$ denotes the fundamental group, which defines a
Novikov ring $\widehat{\mathbb{Z}G}_\xi$, a completion of the group ring. We say that $\omega$
is Morse if it can be represented locally by the differential of a real valued function whose
critical points are nondegenerate. In that case $\omega$ has only finitely many critical points,
each with a well defined index. We write ${\rm ind}\,p$ for the index of $p$.\\[0.2cm]
A vector field $v$ is an $\omega$-gradient if there is a Riemannian metric $g$ such that $\omega_x
(X)=g(X,v(x))$ for every $x\in M$ and $X\in T_xM$. For a critical point $p$ of an $\omega$-gradient
$v$ we denote the unstable, resp. stable, manifold of $p$ by $W^u(p)$, resp. $W^s(p)$.
It is known that $W^u(p)$ is an immersed open disc of
dimension $(n-\mbox{ind }p)$ and $W^s(p)$ one of dimension $\mbox{ind }p$. We say $v$ is
transverse, if all discs $W^s(p)$ and $W^u(q)$ intersect transversely for all critical points
$p,\,q$ of $\omega$. It is well known that in this context one can define a Novikov chain complex
$C_\ast(\omega,v)$ which is in each dimension $i$ a free
$\widehat{\mathbb{Z}G}_\xi$ complex with one generator for every critical point of index $i$.
The boundary homomorphism of $C_\ast(\omega,v)$ is based on the number of trajectories
between critical points of adjacent indices in the universal cover $\tilde{M}$ of $M$.\\[0.2cm]
Suppose such an $\omega$ and $v$ are given, and also a smooth triangulation of $M$. By adjusting
the triangulation if necessary, we can assume that each simplex is transverse to the unstable
manifolds of the critical points of $\omega$. Then there is a natural chain homotopy equivalence
$\varphi(v):\widehat{\mathbb{Z}G}_\xi\otimes_{\mathbb{Z}G}C_\ast^\Delta(\tilde{M})\to
C_\ast(\omega,v)$ given as follows: for a $k$-simplex $\sigma$ we define
\[
\varphi(v)(\sigma)=\sum_{p\in{\rm crit}_k(\omega)}[\sigma:p]\,p
\]
where ${\rm crit}_k(\omega)$ is the set of critical points of $\omega$ having index $k$ and $[
\sigma:p]\in\widehat{\mathbb{Z}G}_\xi$ is the intersection number of a lifting of $\sigma$ to
$\tilde{M}$ with translates of the unstable manifold of a lifting of the critical point $p$.
One asks: what information is contained in the torsion of this equivalence?\\[0.2cm]
For any ring $R$ there is a Dennis trace homomorphism $DT:K_1(R)\to H\!H_1(R)$ from $K$-theory to
Hochschild homology (in all dimensions, but only dimension 1 concerns us here). We use a variant
of $DT$ which we call $\mathfrak{DT}:\overline{W}\to \widehat{H\!H}_1(\mathbb{Z}G)_\xi$.
Here $\overline{W}$ is a subgroup of $K_1(\widehat{\mathbb{Z}G}_\xi)$ containing the torsion of
$\varphi(v)$, and $\widehat{H\!H}_1(\mathbb{Z}G)_\xi$ is a completion of $H\!H_1(\mathbb{Z}G)$ and
related to the Hochschild homology of the Novikov ring by a natural homomorphism $\theta:H\!H_1(
\widehat{\mathbb{Z}G}_\xi)\to\widehat{H\!H}_1(\mathbb{Z}G)_\xi$.\\[0.2cm]
Our main theorem says that if one applies this modified Dennis trace homomorphism to the torsion
of the equivalence $\varphi(v)$, one gets the (topologically important part of the) closed orbit
structure of the flow induced by $v$ in a recognizable form - that of a (noncommutative) zeta
function. In other words the ``Dennis trace'' of the torsion equals the zeta function. Detailed
versions of these sketchy definitions are given in subsequent sections. Here we have just said
enough to state our main theorem.
\begin{theo}\label{itheo}
Let $\omega$ be a Morse form on a smooth connected closed manifold $M^n$. Let $\xi:G\to\mathbb{R}$
be induced by $\omega$ and let $C^\Delta_\ast(\tilde{M})$ be the simplicial $\mathbb{Z}G$ complex
coming from a smooth triangulation of $M$. For every transverse $\omega$-gradient $v$ there is a
natural chain homotopy equivalence $\varphi(v):\widehat{\mathbb{Z}G}_\xi\otimes_{\mathbb{Z}G}
C_\ast^\Delta(\tilde{M})\to C_\ast(\omega,v)$ whose torsion $\tau(\varphi(v))$ lies in
$\overline{W}$ and satisfies
\[\mathfrak{DT}(\tau(\varphi(v)))=\zeta(-v).\]
\end{theo}
Theorem \ref{itheo} is a generalization of \cite[Th.1.1]{schuet}, in fact the whole paper is a
generalization of \cite{schuet} which is closely related to Pajitnov \cite{pajirn,pajitn}. Most
notable is the removal of a cellularity condition on the vector fields $v$ in Theorem \ref{itheo},
a geometric condition due to Pajitnov \cite{pajirn}. Nevertheless this condition is still present
in the proof of Theorem \ref{itheo}. It allows us to identify the torsion of $\varphi(v)$ in a way
that recovers fixed point information which then matches up with closed orbit information. In
order to remove the cellularity condition we show that both the zeta function and the torsion of
$\varphi(v)$ depend continuously on the $\omega$-gradient $v$. The general case then simply follows
from the density of $\omega$-gradients satisfying the cellularity condition among all
$\omega$-gradients.\\[0.2cm]
To define the noncommutative zeta function we use the one parameter fixed point theory of
Geoghegan and Nicas \cite{geonia,geonic} developed for homotopies $F:X\times [a,b]\to X$,
where $X$ is a finite connected CW complex in the sense of classical Nielsen-Wecken fixed point theory.
In the case where $F$ induces the identity on the fundamental group of $X$ they define the one
parameter trace $R(F)$, an algebraically defined element of the Hochschild homology group
$H\!H_1(\mathbb{Z}G)$, where $G$ denotes the fundamental group; if $F$ does not induce the
identity $R(F)$ lies in the first Hochschild homology group with coefficients in a certain
bimodule. This one parameter trace carries information about the fixed points of $F$, i.e. points
$(x,t)\in X\times [a,b]$ with $F(x,t)=x$, and distinguishes between fixed point classes.\\[0.2cm]
Given a vector field on a closed smooth manifold $M$ this theory can be (and is in \cite{geonic})
applied to obtain information about the closed orbit structure of the associated flow by just restricting
the flow to a set $M\times [a,b]$. In the case where this flow only has finitely many closed orbits
in $M\times[a,b]$, the one parameter trace counts these orbits according
to their multiplicity and conjugacy class in $G$.\\[0.2cm]
The noncommutative zeta function of an $\omega$-gradient is now roughly defined as $\zeta(-v)=
\lim_{n\to\infty} R(F_n)$, where $F_n:M\times[0,n]\to M$ is restriction of the flow defined by
$-v$. The vector field has to satisfy a transversality condition so that $\zeta(-v)$ is well
defined.\\[0.2cm]
Commutative zeta functions and their properties
had already been studied by Fried \cite{fried} for homology proper flows. In fact, if the closed
1-form has no critical points, gradient flows are easily seen to be homology
proper, see \cite[\S 2]{fried} for this terminology.\\[0.2cm]
The situation where the closed 1-form is
allowed to have critical points has been studied by Hutchings \cite{hutcth, hutchi}, Hutchings
and Lee \cite{hutlee,hutle2} and Pajitnov \cite{pajirn}, again in a commutative setting.\\[0.2cm]
Noncommutative invariants were studied before that by Geoghegan and Nicas \cite{geonic}
for suspension flows. The case of critical points has only very recently been studied by Pajitnov
\cite{pajitn} for gradient flows
of circle valued maps and by the author \cite{schuet} for gradient flows of closed 1-forms. Instead of a zeta
function, both papers deal with an eta function. In the case where all the closed orbits of $v$ are
nondegenerate, the eta function is defined to be
\[\eta(-v)=\sum\frac{\varepsilon(\gamma)}{m(\gamma)}\{\gamma\}\]
where the sum is taken over the closed orbits $\gamma$ of $-v$, $\varepsilon(\gamma)$ is the Lefschetz sign
of $\gamma$, $m(\gamma)$ its multiplicity and $\{\gamma\}$ the associated conjugacy class in $G$. The
eta function is a well defined object in a quotient of the Novikov ring $
\widehat{\mathbb{R}G}_\xi$, denoted by $\widehat{\mathbb{R}\Gamma}_\xi$. If we replace conjugacy classes
by homology classes, we can take the exponential of the eta function which leads to the commutative zeta
function of Fried \cite{fried}, Hutchings \cite{hutchi,hutcth}, Hutchings and Lee \cite{hutlee,hutle2}
and Pajitnov \cite{pajirn}. But in the noncommutative setting the exponential of the eta function is not
well defined.\\[0.2cm]
The statement of \cite{schuet} and Pajitnov \cite{pajitn} is that, under stronger assumptions on the vector
field $v$, $\eta(-v)=\mathfrak{L}(\tau(\varphi(v)))$, where $\mathfrak{L}$ is a logarithm-like homomorphism
from $\overline{W}$ to $\widehat{\mathbb{R}\Gamma}_\xi$. The connection with Theorem \ref{itheo} is given
by a natural homomorphism $l:\widehat{H\!H}_1(\mathbb{Z}G)_\xi\to\widehat{\mathbb{R}\Gamma}_\xi$ described
in Section \ref{hhofnov} such that $l\circ\mathfrak{DT}=\mathfrak{L}$ and $l(\zeta(-v))=\eta(-v)$.\\[0.2cm]
The contents of this paper are taken from the author's doctoral dissertation written at the State University
of New York at Binghamton under the direction of Ross Geoghegan.
\section{Hochschild homology of group rings}   \label{hhfg}
Let $R$ be a ring and $S$ an $R$-algebra. For an $S-S$ bimodule $M$ we define the \em Hochschild chain complex \em
$(C_\ast(S,M),d)$ by $C_n(S,M)=M\otimes S\otimes\ldots\otimes S$, where the product contains $n$
copies of $S$ and the tensor products are taken over $R$. The boundary operator is given by
\begin{eqnarray*}
d(m\otimes s_1\otimes\ldots\otimes s_n)&=&ms_1\otimes s_2\otimes\ldots\otimes s_n\\
 & &+\sum_{i=1}^{n-1}(-1)^i m\otimes s_1\otimes\ldots\otimes s_is_{i+1}\otimes\ldots\otimes s_n\\
 & &+(-1)^n s_nm\otimes s_1\otimes\ldots\otimes s_{n-1}\end{eqnarray*}
The $n$-th Hochschild homology group of $S$ with coefficients in $M$ is denoted by $H\!H_n(S,M)$.
If $M=S$ and the bimodule structure is given by ordinary multiplication we write $H\!H_\ast(S)$ instead of
$H\!H_\ast(S,M)$. We will mainly be interested in the case where $R=\mathbb{Z}$ and $n=1$. A useful observation is that
$d(x\otimes 1\otimes 1)=x\otimes 1$ and hence classes represented by $x\otimes 1$ are automatically
$0$ in $H\!H_1(S,M)$.\\[0.2cm]
Given an $n\times k$ matrix $A=(A_{ij})$ over $M$ and an $k\times n$ matrix $B=(B_{ij})$ over $S$
we define an $n\times n$ matrix $A\otimes B$ with entries in $M\otimes S$ by setting
$(A\otimes B)_{ij}=\sum_{l=1}^kA_{il}\otimes B_{lj}$. The trace of this matrix, trace $A\otimes B$,
is given by $\sum_{l,m}A_{lm}\otimes B_{ml}$ and is an element of $C_1(S,M)$, it is a cycle if and
only if trace$(AB)=$ trace$(BA)$. Matrices with entries in $C_n(S,M)$ can be defined in a
similar fashion.\\[0.2cm]
Let $G$ be a group and $\phi:G\to G$ an endomorphism. Then we define $(\mathbb{Z}G)^\phi$ to be the
$\mathbb{Z}G-\mathbb{Z}G$ bimodule with underlying abelian group $\mathbb{Z}G$ and multiplication given
by $g\cdot m=\phi(g)m$ and $m\cdot g=mg$ for $m,g\in G$.\\[0.2cm]
We say $g_1$ and $g_2$ in $G$ are \em semiconjugate \em if there exists $g\in G$ with $g_1=\phi(g^{-1})
g_2g$.
We write $\Gamma_\phi$ for the set of semiconjugacy classes and $\gamma(g)$ for the class containing
$g\in G$.\\[0.2cm]
For $\gamma\in\Gamma_\phi$ let $C_n(\mathbb{Z}G,(\mathbb{Z}G)^\phi)_\gamma$ be the subgroup of
$C_n(\mathbb{Z}G,(\mathbb{Z}G)^\phi)$ generated by elements of the form $g_1\otimes g_2\ldots\otimes g_{n+1}$
which satisfy $g_1\cdots g_{n+1}\in\gamma$. Clearly this is a subcomplex and the decomposition
$(\mathbb{Z}G)^\phi\cong\bigoplus_{\gamma\in\Gamma_\phi}\mathbb{Z}\gamma$ as abelian groups gives an isomorphism
of chain complexes
\begin{equation}\label{cgrps}
C_\ast(\mathbb{Z}G,(\mathbb{Z}G)^\phi)\cong\bigoplus_{\gamma\in\Gamma_\phi}C_\ast(\mathbb{Z}G,(\mathbb{Z}G)^\phi)_\gamma.
\end{equation}
Denote the projections by $p_\gamma:C_\ast(\mathbb{Z}G,(\mathbb{Z}G)^\phi)\to C_\ast(\mathbb{Z}G,(\mathbb{Z}G)^\phi)_\gamma$.
It is shown in Geoghegan and Nicas \cite{geonia} that $H_\ast(C_\ast(\mathbb{Z}G,(\mathbb{Z}G)^\phi)_{\gamma(g)})$
is naturally isomorphic to $H_\ast(C(g))$, where $C(g)=\{h\in G|\,g=\phi(h^{-1})gh\}$ denotes the
semicentralizer of $g\in G$.
This gives a natural isomorphism
\begin{equation}\label{hhom}
H\!H_\ast(\mathbb{Z}G,(\mathbb{Z}G)^\phi)\cong\bigoplus_{\gamma(g)\in\Gamma_\phi}H_\ast(C(g)).
\end{equation}
In particular we have $H\!H_0(\mathbb{Z}G,(\mathbb{Z}G)^\phi)\cong\mathbb{Z}\Gamma_\phi$ and
$H\!H_1(\mathbb{Z}G,(\mathbb{Z}G)^\phi)\cong\bigoplus_{\gamma(g)\in\Gamma_\phi}C(g)_{ab}$, the direct
sum of the abelianizations of the semicentralizers.
\section{Review of parametrized fixed point theory}\label{opfpt}
In this section we recall the parametrized fixed point theory of Geoghegan and Nicas \cite{geonia}.
Let $X$ be a finite connected CW complex, $v\in X$ a basepoint and $F:X\times I^r\to X$ a cellular
map, where $r\geq 0$ and $I^r$ is the product of $r$ copies of the unit interval with the usual CW
structure. We set $G=\pi_1(X,v)$. Then the choice of a basepath $\tau$ from $v$ to $F(v,0,\ldots\hspace{-1pt},0)$
induces an endomorphism $\phi:G\to G$ given by $\phi([\omega])=[\tau\ast F_0\circ\omega\ast\tau^{-1}]$.
Let $\tilde{X}$ be the universal covering space of $X$ and $\tilde{v}$ a lift of $v$.
Let $\tilde{\tau}$ be the lift of $\tau$ which starts at $\tilde{v}$ and let $\tilde{F}$ be the lift of
$F$ mapping $(\tilde{v},0,\ldots\hspace{-1pt},0)$ to $\tilde{\tau}(1)$.
By choosing an oriented
lift $\tilde{e}$ for every cell $e$ in $X$ we get a basis of the free left $\mathbb{Z}G$ complex $C_\ast
(\tilde{X})$. The action is given by $[\omega]\tilde{e}=h_{[\omega]\ast}(\tilde{e})$, where
$h_{[\omega]}$ is the covering transformation that sends $\tilde{v}$ to $\tilde{\omega}(1)$, where
$\tilde{\omega}$ is a lift of $\omega$ with $\tilde{\omega}(0)=\tilde{v}$.
\begin{rema}\em In \cite{geonia}, Geoghegan and Nicas consider $C_\ast(\tilde{X})$ as a right complex.
This leads to differences between their exposition and ours in that group elements have to be
replaced here by their inverse, semiconjugation has the $\phi$ on the left, etc. This has no impact
on the main theorems in \cite{geonia,geonic} other than sign differences.
\end{rema}
We now define $\tilde{D}^F_k:C_k(\tilde{X})\to C_{k+r}(\tilde{X})$ by $\tilde{e}\tilde{D}_k^F=
(-1)^{(k+r)r}\tilde{F}_{k+r}(\tilde{e}\times I^r)$. Since we consider $C_\ast(\tilde{X})$ as a left
module, we write $\tilde{D}_k^F$ on the right. A standard computation gives \\
$(g\tilde{e})\tilde{D}_k^F=\phi(g)(\tilde{e}\tilde{D}_k^F)$. We want to examine the behavior with
the boundary homomorphism of $C_\ast(\tilde{X})$. So assume that $r\geq 1$. Then
\begin{eqnarray*}
\tilde{e}\tilde{\partial}_k\tilde{D}^F_{k-1}&=&(-1)^{(k-1+r)r}\tilde{F}_{k+r-1}(\tilde{\partial}_k\tilde{e}\times I^r)\\
\tilde{e}\tilde{D}^F_k\tilde{\partial}_{k+r}&=&(-1)^{(k+r)r}\tilde{\partial}_{k+r}\tilde{F}_{k+r}(\tilde{e}\times I^r)\\
&=&(-1)^{(k+r)r}\tilde{F}_{k+r-1}(\tilde{\partial}_k\tilde{e}\times I^r)+(-1)^{(k+r)r+k}\tilde{F}_{k+r-1}(\tilde{e}\times \partial I^r)\\
&=&(-1)^{(k+r)r}\tilde{F}_{k+r-1}(\tilde{\partial}_k\tilde{e}\times I^r)-\sum_{j=1}^{2^r}(-1)^{\sigma(j)}\tilde{D}^{F_j}_k(\tilde{e}).
\end{eqnarray*}
Here $F_j:X\times I^{r-1}\to X$ is obtained from $F$ by restricting to a side of $\partial I^r$. The
sign $(-1)^{\sigma(j)}$ depends on the orientation of that side in $\partial I^r$.
Hence
\[\tilde{D}^F_k\tilde{\partial}_{k+r}+(-1)^{r+1}\tilde{\partial}_k\tilde{D}^F_{k-1}=\sum_{j=1}^{2^r}
(-1)^{\sigma(j)}\tilde{D}^{F_j}_k.\]
We define endomorphisms of $C_\ast(\tilde{X})$ by $\tilde{D}_\ast^F=\bigoplus_k(-1)^{k+r}\tilde{D}^F_k$,
$\tilde{D}_\ast^{F_j}=\bigoplus_k(-1)^{k+r-1}\tilde{D}^{F_j}_k$ and $\tilde{\partial}_\ast=\bigoplus_k
\tilde{\partial}_k$ and denote the matrices with the same letter. The alternation of signs leads
to the following matrix equality:
\begin{equation}\label{matrix}
\tilde{D}^F_\ast\tilde{\partial}_\ast+(-1)^r\phi(\tilde{\partial}_\ast)\tilde{D}^F_\ast=-
\sum_{j=1}^{2^r}(-1)^{\sigma(j)}\tilde{D}^{F_j}_\ast.
\end{equation}
For $r=1$ this implies
\begin{equation} \label{bound}
d(\mbox{trace}(\tilde{D}^F_\ast\otimes\tilde{\partial}_\ast))=\mbox{trace}(\tilde{D}^F_\ast
\tilde{\partial}_\ast)-\mbox{trace}(\phi(\tilde{\partial}_\ast)\tilde{D}^F_\ast)=
\mbox{trace }\tilde{D}^{F_0}_\ast-\mbox{trace }\tilde{D}^{F_1}_\ast.
\end{equation}
These traces contain information about the fixed points of the respective maps. To clarify this
let Fix$(F)=\{(x,t_1,\ldots\hspace{-1pt},t_r)\in X\times I^r|\, F(x,t_1,\ldots\hspace{-1pt},t_r)=x\}$. We define an equivalence
relation $\sim$ on Fix$(F)$ by saying that the fixed points $a$ and $b$ are equivalent
if there exists a path $\nu$ in $X\times I^r$ from $a$ to $b$ such that the loop $(p\circ\nu)\ast(F\circ\nu)^{-1}$
is homotopically trivial; here $p:X\times I^r\to X$ is projection. Since $X\times I^r$ is compact and
locally contractible there are only finitely many fixed point classes.
There is an injective function $\Phi:\mbox{Fix}(F)/\!\sim\,\to\Gamma_\phi$ by mapping $[a]$ to $\gamma([\tau\ast
(F\circ\mu)\ast(p\circ\mu)^{-1}])$, where $\mu$ is a path from $(v,0,\ldots\hspace{-1pt},0)$ to $a$. It is easy to see
that this is well defined, but it depends on the choice of the base path $\tau$. \\[0.2cm]
Let us return to the algebraic traces. First we look at a cellular map $f:X\to X$. Geoghegan and Nicas
\cite{geonia} prove the following
\begin{theo}\cite[Th.2.6]{geonia}
Let $f:X\to X$ be a cellular self map of a CW complex $X$. Let $e$ be a $q$-cell of $X$ and let the
corresponding diagonal entry in the $\mathbb{Z}G$ matrix of $\tilde{f}_q:C_q(\tilde{X})\to C_q(\tilde{X})$
be $d(\tilde{e})=\sum gm_g$, where $m_g\in\mathbb{Z}$. For each $m_g\not=0$, $e$ contains a fixed
point $x_g$ such that $\tilde{f}(\tilde{x}_g)=g\tilde{x}_g$, where $\tilde{x}_g$ is the lift of
$x_g$ to $\tilde{e}$.
\end{theo}
In particular, nonzero terms in $p_\gamma($trace $\tilde{D}^f_\ast)\in\mathbb{Z}\gamma\subset(\mathbb{Z}G)^\phi$
detect fixed points $x$ with $\Phi[x]=\gamma$.
In other words if $f$ has no fixed points corresponding to $\gamma$, then $p_\gamma($trace $\tilde{D}^f_\ast)=0$.
This leads to the classical Nielsen-Wecken fixed point theory; see Geoghegan and Nicas \cite{geonia}
for details.\\[0.2cm]
To get a similar result for one parameter fixed point theory we need the following theorem of
Geoghegan and Nicas \cite{geonia}.
\begin{theo}\cite[Th.2.12]{geonia}
Let $F:X\times I\to X$ be a cellular map where $X$ is a well faced CW complex. Let $e^{q-1}\subset e^q$
be cells of $X$ of the indicated dimensions and let the corresponding entry in the $\mathbb{Z}G$ matrix
of $\tilde{F}_q:C_q(\tilde{X}\times I)\to C_q(\tilde{X})$ be
$d(\tilde{e}^{q-1},\tilde{e}^q)=\sum gm_g$ where each $m_g\in\mathbb{Z}$. For each $g\in G$ with $m_g\not=0$
and $u\in G$ with $h_u(\tilde{e}^{q-1})\subset\tilde{e}^q$, $e^{q-1}\times I$ contains a
fixed point $(x_{g,u},t_{g,u})$ such that $\tilde{F}(\tilde{x}_{g,u},t_{g,u})=gu\tilde{x}_{g,u}$,
where $(\tilde{x}_{g,u},t_{g,u})$ is the lift of $(x_{g,u},t_{g,u})$ to $\tilde{e}^{q-1}\times I$.
\end{theo}
The condition ``well faced" is merely a technicality which is satisfied for example by regular CW
complexes. Again we get that nonzero terms in $p_\gamma($trace$(\tilde{D}^F_\ast\otimes\tilde{\partial}_\ast))
\in C_1(\mathbb{Z}G,(\mathbb{Z}G)^\phi)_\gamma$ detect fixed points $(x,t)$ with $\Phi[x,t]=\gamma$.\\[0.2cm]
Looking back at (\ref{bound}) we now see that $p_\gamma($trace$(\tilde{D}^F_\ast\otimes\tilde{\partial}_\ast))$
is a cycle if $F_1$ and $F_0$ have no fixed points associated to $\gamma$. Let $S\subset\Gamma_\phi$, then
we define $p^\dag_S:C_\ast(\mathbb{Z}G,(\mathbb{Z}G)^\phi)\to C_\ast(\mathbb{Z}G,(\mathbb{Z}G)^\phi)$ to
be the composition
\[ C_\ast(\mathbb{Z}G,(\mathbb{Z}G)^\phi)\stackrel{\bigoplus\limits_{\gamma\in\Gamma_\phi-S}\,p_\gamma}{\longrightarrow}
\bigoplus_{\gamma\in\Gamma_\phi-S}C_\ast(\mathbb{Z}G,(\mathbb{Z}G)^\phi)_\gamma
\stackrel{i}{\hookrightarrow}\bigoplus_{\gamma\in\Gamma_\phi}C_\ast(\mathbb{Z}G,(\mathbb{Z}G)^\phi)_\gamma
=C_\ast(\mathbb{Z}G,(\mathbb{Z}G)^\phi).\]
\begin{defi}\em Let $F:X\times I\to X$ be a cellular homotopy, where $X$ is a finite connected CW
complex. Then the \em one parameter trace \em of $F$, denoted by $R(F)\in H\!H_1(\mathbb{Z}G,(\mathbb{Z}G)^\phi)$,
is the homology class of $p^\dag_S($trace$(\tilde{D}^F_\ast\otimes\tilde{\partial}_\ast))$, where $S\subset\Gamma_\phi$
is the set of semiconjugacy classes associated to fixed points of $F_0$ and $F_1$.
\end{defi}
By the remarks above we get that if $p_{\gamma(g)\ast}R(F)\in C(g)_{ab}$ is nonzero, then $F$ contains
a fixed point associated to $\gamma(g)$.
It is shown in Geoghegan and Nicas \cite[Prop.4.1]{geonia} that $R(F)$ is independent of the
orientation and the choice of lifts of cells to $\tilde{X}$.
The following is an easy example, but it shows how the signs behave compared to Geoghegan and Nicas
\cite{geonia,geonic}.
\begin{exam}\em Let $\Phi:S^1\times\mathbb{R}\to S^1$ be defined by $\Phi(e^{2\pi i\theta},t)=e^{2\pi i(\theta
+t)}$ and let \\$F_n:S^1\times[0,n]\to S^1$ be given by $F_n=\Phi|_{S^1\times[0,n]}$. The basepath is chosen
to be constant. We can put a cell structure on $S^1$ with two cells and lift it to a cell structure
of $\mathbb{R}=\tilde{S^1}$. Choose a $0$-cell $\tilde{e}^0=0\in\mathbb{R}$ and a 1-cell $\tilde{e}^1
=[0,1]$. Let $t$ be the generator of $\pi_1(S^1,1)=G$ that satisfies $h_t(x)=x+1$.
Then $\tilde{e}^0\tilde{D}_0^F=-(1+t+\ldots+t^{n-1})\tilde{e}^1$ and $\tilde{e}^1\tilde{\partial}_1=(t-1)\tilde{e}^0$.
Hence $\mbox{trace}(\tilde{D}^F_\ast\otimes\tilde{\partial}_\ast)=(1+t+\ldots+t^{n-1})\otimes(t-1)$ which is
homologous to $1\otimes t + t\otimes t+\ldots+t^{n-1}\otimes t$. The last summand corresponds to a fixed point of $F(\cdot,n)$,
so $R(F_n)$ is represented by $\sum_{k=1}^{n-1}t^{k-1}\otimes t$. Notice that $H_1(C_\ast(\mathbb{Z}G,
\mathbb{Z}G)_{\gamma(t^k)})
\simeq H_1(C(t^k))\simeq\mathbb{Z}$ and $t^{k-1}\otimes t$ is a generator.
\label{examplos}\end{exam}
To examine the behavior of $R(F)$
within its homotopy class we look at chains modulo boundaries. Define $C/B=C_1(\mathbb{Z}G,(\mathbb{Z}G)^\phi)
/d_2(C_2(\mathbb{Z}G,(\mathbb{Z}G)^\phi))$ and let $CR(F)$ be the image of trace$(\tilde{D}^F_\ast
\otimes\tilde{\partial}_\ast)$ in $C/B$. Given a cellular homotopy $\Lambda:X\times I\times I\to X$
define \\$F^i(x,t)=\Lambda(x,t,i)$ and $U^i(x,t)=\Lambda(x,i,t)$ for $i=0,1$. To ensure that all maps
induce the same homomorphism $\phi$, $F^0$ and $U^0$ use the basepath $\tau$ while $F^1$ uses $\tau\ast\sigma$
and $U^1$ uses $\tau\ast\mu$, where $\sigma(t)=\Lambda(v,0,t)$ and $\mu(t)=\Lambda(v,t,0)$. Then
we get the following
\begin{prop}\cite[Prop.4.2]{geonia} $CR(F^0)-CR(F^1)=CR(U^0)-CR(U^1)$. \label{hprop}
\end{prop}
\begin{proof}
Use (\ref{matrix}) with $r=2$. This gives
\begin{eqnarray*}d(\mbox{trace}(\tilde{D}^\Lambda_\ast\otimes\tilde{\partial}_\ast\otimes\tilde{\partial}_\ast))&
=&\mbox{trace}((\tilde{D}^\Lambda_\ast\tilde{\partial}_\ast+
\phi(\tilde{\partial}_\ast)\tilde{D}^\Lambda_\ast)\otimes\tilde{\partial}_\ast)\\
&=&\mbox{trace}(\tilde{D}^{U^0}_\ast\otimes\tilde{\partial}_\ast)-
 \mbox{trace}(\tilde{D}^{U^1}_\ast\otimes\tilde{\partial}_\ast)+ \\
&&\mbox{trace}(\tilde{D}^{F^1}_\ast\otimes\tilde{\partial}_\ast)-
 \mbox{trace}(\tilde{D}^{F^0}_\ast\otimes\tilde{\partial}_\ast).
\end{eqnarray*}
Passing to $C/B$ gives the result.
\end{proof}
In particular $R(F)$ only depends on the homotopy class of $F$ relative to the ends. Another property
of $R(F)$ is combinatorial invariance; see \cite[Cor.4.6]{geonia} and \cite[Th.4.6]{genisc}, which
is used in \cite[\S 4B]{geonia} to extend the definition of $R(F)$ to continuous homotopies on compact polyhedra.
To do this one just has to use a fine enough triangulation and a simplicial approximation of $F$.
This $R(F)$ has the same properties for detecting fixed points and homotopy invariance.
\section{Hochschild homology of Novikov rings}\label{hhofnov}
Let $G$ be a group and $\xi:G\to\mathbb{R}$ be a homomorphism. We denote by $\widehat{\widehat{\mathbb{Z}G}}$
the abelian group of all functions $G\to\mathbb{Z}$. For $\lambda\in\widehat{\widehat{\mathbb{Z}G}}$
let supp $\lambda=\{g\in G\,|\,\lambda(g)\not=0\}$. Then we define
\[\widehat{\mathbb{Z}G}_\xi=\{\lambda\in\widehat{\widehat{\mathbb{Z}G}}\,|\,\forall r\in\mathbb{R}
\hspace{0.4cm}\#\,\mbox{supp }\lambda\cap\xi^{-1}([r,\infty))<\infty\}\]
For $\lambda_1,\lambda_2\in\widehat{\mathbb{Z}G}_\xi$ we set $(\lambda_1\cdot\lambda_2)(g)=\sum
\limits_{h_1,h_2\in G\atop h_1h_2=g}\lambda_1(h_1)\lambda_2(h_2)$, then $\lambda_1\cdot\lambda_2$
is a well defined element of $\widehat{\mathbb{Z}G}_\xi$ and turns $\widehat{\mathbb{Z}G}_\xi$ into
a ring, the \em Novikov ring\em. It contains the usual group ring $\mathbb{Z}G$ as a subring and
we have $\mathbb{Z}G=\widehat{\mathbb{Z}G}_\xi$ if and only if $\xi$ is the zero homomorphism.
We can define Novikov rings for $\mathbb{Q}$ or $\mathbb{R}$ by simply replacing
$\mathbb{Z}$ in the above definitions by one of these.
\begin{defi}\label{dnorm}\em The \em norm \em of $\lambda\in\widehat{\mathbb{Z}G}_\xi$ is defined to be
\[\|\lambda\|=\|\lambda\|_\xi=\inf\{t\in(0,\infty)|\mbox{ supp }\lambda\subset\xi^{-1}((-\infty,\log t])\}.\]
\end{defi}
It has the following nice properties:
\begin{enumerate}
\item $\|\lambda\|\geq 0$ and $\|\lambda\|=0$ if and only if $\lambda=0$.
\item $\|\lambda\|=\|-\lambda\|$.
\item $\|\lambda+\mu\|\leq\max\{\|\lambda\|,\|\mu\|\}$.
\item $\|\lambda\cdot\mu\|\leq\|\lambda\|\cdot\|\mu\|$.
\end{enumerate}
This norm can be used to define a complete metric on $\widehat{\mathbb{Z}G}_\xi$, the topology induced by this metric
is the Krull topology, compare Eisenbud \cite{eisenb}.\\[0.2cm]
If $N$ is a normal subgroup of $G$ that is contained in $\ker\xi$ we get a well defined homomorphism
$\bar{\xi}:G/N\to\mathbb{R}$ and a well defined ring epimorphism $\varepsilon:\widehat{\mathbb{Z}G}_\xi
\to\widehat{\mathbb{Z}G/N}_{\bar{\xi}}$ given by \\$\varepsilon(\lambda)(gN)=\sum\limits_{n\in N}\lambda(gn)$.\\[0.2cm]
Now let $\Gamma=\Gamma_{1}$ be the set of conjugacy classes of $G$. Again the homomorphism $\xi$
induces a well defined function $\Gamma\to\mathbb{R}$ which we also denote by $\xi$. In analogy with above
we define $\widehat{\mathbb{Z}\Gamma}_\xi$, but since there is no well defined multiplication in
$\Gamma$, this object is just an abelian group. Again there is an epimorphism $\varepsilon:
\widehat{\mathbb{Z}G}_\xi\to\widehat{\mathbb{Z}\Gamma}_\xi$ of abelian groups. We can think of
$\widehat{\mathbb{Z}\Gamma}_\xi$ as lying between $\widehat{\mathbb{Z}G}_\xi$ and
$\widehat{\mathbb{Z}H_1(G)}_{\bar{\xi}}$.\\[0.2cm]
We would like to get a result for $H\!H_\ast(\widehat{\mathbb{Z}G}_\xi)$ similar to (\ref{hhom}).
We will just be interested in the case where $\phi$ is the identity. Before we go into detail let
us have an informal discussion. Elements of $\mathbb{Z}G$ and $\widehat{\mathbb{Z}G}_\xi$
can also be described as formal linear combinations. So a typical $n$-chain looks like
\[\sum_{g_1\in G}n_{g_1}g_1\otimes\ldots\otimes\sum_{g_{n+1}\in G}n_{g_{n+1}}g_{n+1}.\]
If the elements are taken from $\mathbb{Z}G$ we can write this as
\begin{equation}\label{schwups}
\sum_{g_1,\ldots\hspace{-1pt},g_{n+1}\in G}n_{g_1}\cdots n_{g_{n+1}}\,g_1\otimes \ldots\otimes g_{n+1}
\end{equation}
which is just a finite sum. But if we think of the elements as being taken from $\widehat{\mathbb{Z}G}
_\xi$, (\ref{schwups}) would be an infinite sum of tensors which does not give a well defined
element of $C_n(\widehat{\mathbb{Z}G}_\xi,\widehat{\mathbb{Z}G}_\xi)$. In other words, the process
of breaking down an $n$-chain into the form (\ref{schwups}) might not give an $n$-chain. On the
other hand given a conjugacy class $\gamma\in\Gamma$ there are only finitely many nonzero
summands in (\ref{schwups}) that satisfy $g_1\cdots g_{n+1}\in\gamma$. So we are going to define
a chain complex $C_\ast$ based on conjugacy classes in which (\ref{schwups}) makes sense together
with a chain map $C_n(\widehat{\mathbb{Z}G}_\xi,\widehat{\mathbb{Z}G}_\xi)\to C_\ast$.\\[0.2cm]
More precisely for $\gamma\in\Gamma$ we define
a chain map $\theta_\gamma:C_\ast(\widehat{\mathbb{Z}G}_\xi,\widehat{\mathbb{Z}G}_\xi)\to C_\ast
(\mathbb{Z}G,\mathbb{Z}G)_\gamma$. To do this let $\lambda_1\otimes\lambda_2\ldots\otimes\lambda_{n+1}
\in C_n(\widehat{\mathbb{Z}G}_\xi,\widehat{\mathbb{Z}G}_\xi)$. We can assume that all $\lambda_i$ are
nonzero. Then let $T=\log\min\{\|\lambda_i\||\,i=1,\ldots\hspace{-1pt},n+1\}$ and choose an $M<0$ such that
\\$\|\lambda_1\|\cdots\|\lambda_{n+1}\|\leq\exp(-M)$. Then for $i=1,\ldots\hspace{-1pt},n+1$ we define $\bar{\lambda}_i
\in\mathbb{Z}G$ by
\[\bar{\lambda}_i(h)=\left\{\begin{array}{rl}0&\mbox{if }\xi(h)<M+T+\xi(\gamma)\\
\lambda_i(h)&\mbox{otherwise}\end{array}\right.\]
and
\[\theta_\gamma(\lambda_1\otimes\ldots\otimes\lambda_{n+1})=p_\gamma(\bar{\lambda}_1\otimes
\ldots\otimes\bar{\lambda}_{n+1}).\]
It is to be shown that this is independent of $M$, so let $\bar{\bar{\lambda}}_i$ be defined using
an $M'$ with $M'<M$. It suffices to show that $p_\gamma(\bar{\bar{\lambda}}_1\otimes\ldots\otimes
\bar{\bar{\lambda}}_{i-1}\otimes\bar{\bar{\lambda}}_i-\bar{\lambda}_i\otimes\bar{\lambda}_{i+1}\otimes
\ldots\otimes\bar{\lambda}_{n+1})=0$. Let $g_i\in\mbox{supp}(\bar{\bar{\lambda}}_i-\bar{\lambda}_i)$
and assume that there are $g_j\in\mbox{supp}\,\bar{\bar{\lambda}}_j$ for $j\not=i$ with $g_1\cdots
g_{n+1}\in\gamma$. Then
\[\xi(g_i)<M+T+\xi(\gamma)=M+T+\xi(g_1)+\ldots +\xi(g_{n+1}),\]
so
\[\sum_{j\not= i}\xi(g_j)>-M-T,\]
hence
\[\exp(\sum_{j\not= i}\xi(g_j))>\exp(-M)\cdot\exp(-T)\geq\exp(-M)\cdot\|\lambda_i\|^{-1}\]
and so $\|\lambda_1\|\cdots\|\lambda_{n+1}\|>\exp(-M)$, a contradiction, hence the $g_j$ for $j\not= i$
with the described properties cannot exist. Therefore $p_\gamma(\bar{\bar{\lambda}}_1\otimes\ldots\otimes
\bar{\bar{\lambda}}_{i-1}\otimes\bar{\bar{\lambda}}_i-\bar{\lambda}_i\otimes\bar{\lambda}_{i+1}\otimes
\ldots\otimes\bar{\lambda}_{n+1})=0$.\\[0.2cm]
Next we show that $\theta_\gamma$ commutes with the boundary.
We have
\begin{eqnarray*}
d(\lambda_1\otimes\ldots\otimes\lambda_{n+1})&=&\lambda_1\lambda_2\otimes\ldots\otimes \lambda_{n+1}\\
 & &+\sum_{i=1}^{n-1}(-1)^i \lambda_1\otimes\ldots\otimes \lambda_i\lambda_{i+1}\otimes\ldots\otimes \lambda_{n+1}\\
 & &+(-1)^n \lambda_{n+1}\lambda_1\otimes\ldots\otimes \lambda_n\end{eqnarray*}
In forming `` $\bar{ }$ '', use $T=\log(\min\{\|\lambda_i\||\,i=1,\ldots\hspace{-1pt},n+1\}\cup\{\|\lambda_i\|\cdot\|
\lambda_j\||\,i,j=1,\ldots\hspace{-1pt},n+1\})$ and $M$ so small that the same $M+T+\xi(\gamma)$ can be used for every
$\bar{\lambda}_i,\overline{\lambda_i\lambda_{i+1}}$. We have to show that $p_\gamma(\bar{\lambda}_1
\otimes\ldots\otimes\overline{\lambda_i\lambda_{i+1}}-\bar{\lambda}_i\bar{\lambda}_{i+1}\otimes
\ldots\otimes\bar{\lambda}_{n+1})=0$. We have
\[\overline{\lambda_i\lambda_{i+1}}(h)=\left\{\begin{array}{rl}0&\mbox{if }\xi(h)<M+T+\xi(\gamma)\\
\lambda_i\lambda_{i+1}(h)&\mbox{otherwise}\end{array}\right.\]
On the other hand
\[\bar{\lambda}_i\bar{\lambda}_{i+1}(h)=\sum_{h_ih_{i+1}=h}\bar{\lambda}_i(h_i)\bar{\lambda}_{i+1}(h_{i+1})
=\sum_{h_ih_{i+1}=h\atop \xi(h_k)\geq M+T+\xi(\gamma)}\lambda_i(h_i)\lambda_{i+1}(h_{i+1}).\]
Let $h\in\mbox{supp}(\overline{\lambda_i\lambda_{i+1}}-\bar{\lambda}_i\bar{\lambda}_{i+1})$. If
$\xi(h)<M+T+\xi(\gamma)$, then there exist $h_i\in\mbox{supp}\,\lambda_i$, $h_{i+1}\in\mbox{supp}\,\lambda_{i+1}$
with $h_ih_{i+1}=h$. The existence of $h_j\in\mbox{supp}\,\lambda_j$ for $j\not=i,i+1$ with $h_1\cdots h_{n+1}\in\gamma$
leads to a contradiction as above. If $\xi(h)\geq M+T+\xi(\gamma)$, then assume without loss of
generality that there is $h_i\in\mbox{supp}\,\lambda_i$, $h_{i+1}\in\mbox{supp}\,\lambda_{i+1}$ with
$\xi(h_i)<M+T+\xi(\gamma)$ and $\xi(h_{i+1})\geq M+T+\xi(\gamma)$, but as before no $h_j\in\mbox{supp}\,
\lambda_j$ for $j\not=i,i+1$with $h_1\cdots h_{n+1}\in\gamma$ can exist.
Therefore $p_\gamma(\bar{\lambda}_1
\otimes\ldots\otimes\overline{\lambda_i\lambda_{i+1}}-\bar{\lambda}_i\bar{\lambda}_{i+1}\otimes
\ldots\otimes\bar{\lambda}_{n+1})=0$ and $\theta_\gamma$ is a chain homomorphism.\\[0.2cm]
In analogy with (\ref{cgrps}) we define
\[C_\ast(\mathbb{Z}G)_\xi=\{(c_\gamma)\in\prod_{\gamma\in\Gamma}C_\ast(\mathbb{Z}G,\mathbb{Z}G)_\gamma|
\,\forall r\in\mathbb{R}
\hspace{0.4cm}\#\,\{ c_\gamma\not=0\,|\,\xi(\gamma)\geq r\}<\infty\}\]
and $\widehat{H\!H}_\ast(\mathbb{Z}G)_\xi=H_\ast(C_\ast(\mathbb{Z}G)_\xi)$. Notice that
\[\bigoplus_{\gamma\in\Gamma}H_\ast(C_\ast(\mathbb{Z}G,\mathbb{Z}G)_\gamma)\subset\widehat{H\!H}_\ast(\mathbb{Z}G)_\xi\subset\prod
_{\gamma\in\Gamma}H_\ast(C_\ast(\mathbb{Z}G,\mathbb{Z}G)_\gamma),\]
which follows from Section \ref{hhfg}. Furthermore it is easy to see that \\$(\theta_{\gamma\ast})
_{\gamma\in\Gamma}:H\!H_\ast(\widehat{\mathbb{Z}G}_\xi)\to\prod_{\gamma\in\Gamma}H_\ast(C_\ast(\mathbb{Z}G,\mathbb{Z}G)_\gamma)$
factors through $\widehat{H\!H}_\ast(\mathbb{Z}G)_\xi$; we denote this corestriction by $\theta$.
Therefore we have the commutative diagram
\[
\begin{array}{ccc}
H\!H_\ast(\mathbb{Z}G)&\stackrel{\simeq}{\longrightarrow}&\bigoplus\limits_{\gamma\in\Gamma}H_\ast(C_\ast(\mathbb{Z}G,\mathbb{Z}G)_\gamma)\\[0.4cm]
\Big\downarrow& &\Big\downarrow\\[0.3cm]
H\!H_\ast(\widehat{\mathbb{Z}G}_\xi)&\stackrel{\theta}{\longrightarrow}&\widehat{H\!H}_\ast(\mathbb{Z}G)_\xi
\end{array}
\]
By analogy with \cite[\S 3]{schuet} we can define a homomorphism $L:\prod\limits_{\gamma\in\Gamma}C_1(\mathbb{Z}G,\mathbb{Z}G)_\gamma
\to\widehat{\widehat{\mathbb{R}\Gamma}}$
by
\[L((g_1\otimes g_2)_{g_1g_2\in\gamma})(\gamma)=\left\{\begin{array}{rl}
\dfrac{\xi(g_2)}{\xi(\gamma)}&\mbox{if }\xi(\gamma)<0\\
0&\mbox{otherwise}\end{array}\right.\]
which induces a homomorphism on homology that restricts to a homomorphism 
\[l:\widehat{H\!H}_1(\mathbb{Z}G)_\xi\to\widehat{\mathbb{R}\Gamma}_\xi.\]
In \cite[\S 3]{schuet}, the bimodule in the Hochschild complex is on the
right, so if the homomorphism $\mu$ in \cite{schuet} is defined by switching
the bimodule to the left, we see that this homomorphism factors as
$\mu=l\circ\theta$. \begin{rema}\label{remarquis}\em In Section \ref{hhfg} we
saw that $H\!H_0(\mathbb{Z}G)\cong\mathbb{Z}\Gamma$, a result which can easily
be derived directly, in particular we also have $H\!H_0(\mathbb{R}G)\cong
\mathbb{R}\Gamma$. We can also define $C_n(\mathbb{R}G,\mathbb{R}G)_\gamma$ to
be the subgroup of $C_n(\mathbb{R}G,\mathbb{R}G)$ generated by
$r_1g_1\otimes\ldots\otimes r_{n+1}g_{n+1}$ with $r_1,
\ldots,\hspace{-1pt},r_{n+1}\in\mathbb{R}$ and $g_1\cdots g_{n+1}\in\gamma$.
In analogy with above we get the complex $C_\ast(\mathbb{R}G)_\xi$, which
contains $C_\ast(\mathbb{R}G,\mathbb{R}G)$ as a subcomplex. Denoting the
resulting homology by $\widehat{H\!H}_\ast(\mathbb{R}G)_\xi$, a completion of
$H\!H_\ast(\mathbb{R}G)$, we get
$\widehat{\mathbb{R}\Gamma}_\xi\cong\widehat{H\!H}_0(\mathbb{R}G)_\xi$ and the
homomorphism $l$ is reminiscent of the homomorphism $\hat{P}_+$ in Geoghegan
and Nicas \cite[\S 5]{geonic}. \end{rema}
\section{Gradient flows of closed 1-forms}\label{gradflws}
Given a closed 1-form $\omega$ on a closed manifold $M$ we obtain a
homomorphism $\bar{\xi}:H_1(M)\to \mathbb{R}$ by
$\bar{\xi}[\alpha]=\int_\alpha\omega$ which induces a homomorphism
$\xi:G\to\mathbb{R}$ , where $G=\pi_1(M,v)$ for some basepoint $v\in M$. Since
$G$ is finitely generated, the image of $\xi$ is a finitely generated subgroup
of $\mathbb{R}$, hence isomorphic to $\mathbb{Z}^k$ for some integer $k$. If
$k=1$ $\omega$ is said to be \em rational\em, if $k>1$ it is \em
irrational\em.\\[0.2cm] We will call a closed $1$-form a \em Morse form \em if
$\omega$ is locally represented by the differential of real valued functions
whose critical points are nondegenerate. So if $\omega$ is a Morse form, then
$\omega$ has only finitely many critical points and every critical point has a
well defined index. If $p$ is a critical point, we denote its index by
$\mbox{ind }p$.
\begin{defi}\em Let $\omega$ be a closed $1$-form. A vector
field $v$ is called an \em $\omega$-gradient\em, if there exists a Riemannian
metric $g$ such that $\omega_x(X)=g(X,v(x))$ for every $x\in M$ and $X\in T_x
M$.
\end{defi}
For a critical point $p$ of an $\omega$-gradient $v$ we denote the
unstable, resp. stable, manifold of $p$ by $W^u(p)$, resp. $W^s(p)$. So if
$\Phi:M\times\mathbb{R}\to M$ denotes the flow of $v$, then \\$W^u(p)=\{x\in
M|\Phi(x,t)\to p$ for $t\to-\infty\}$ and $W^s(p)=\{x\in M|\Phi(x,t)\to p$ for
$t\to\infty\}$. It is known that $W^u(p)$ is an immersed open disc of
dimension $(n-\mbox{ind }p)$ and $W^s(p)$ one of dimension $\mbox{ind }p$, see e.g.\ Abraham and
Robbin \cite[\S 27]{abrah}.
The next Lemma allows us to forget about the Riemannian metric
and will be useful in using vector fields as gradients of different Morse
forms.
\begin{lem}\label{smalo}\cite[Lm.2.3]{schuet}
Let $\omega$ be a Morse
form and $v$ a vector field. Then $v$ is an $\omega$-gradient if and only if
\begin{enumerate} \item For every critical point $p$ of $\omega$ there exists
a neighborhood $U_p$ of $p$ and a Riemannian metric $g$ on $U_p$ such that
$\omega_x(X)=g(X,v(x))$ for every $x\in U_p$ and $X\in T_x U_p$. \item If
$\omega_x\not=0$, then $\omega_x(v(x))>0$. \end{enumerate}
\end{lem}
Condition 2.\ will sometimes be all we need so we define a smooth vector field $v$ to be a
\em weak $\omega$-gradient \em if $\omega_x=0$ implies $v(x)=0$ and $\omega_x\not=0$ implies
$\omega_x(v(x))>0$. The stable and unstable manifolds still exist as sets but they might not
have as nice properties as in the case of $\omega$-gradients.
\begin{defi}\label{transdf}
\em Let $v$ be an $\omega$-gradient.
\begin{enumerate}
\item We say $v$ is \em transverse\em, if all discs $W^s(p)$ and $W^u(q)$ intersect
transversely for all critical points $p,q$ of $\omega$.
\item We say $v$ is \em almost transverse\em, if for critical points $p,q$ with
$\mbox{ind }p\leq\mbox{ind } q$ the discs $W^s(p)$ and $W^u(q)$ intersect transversely.
\end{enumerate}
\end{defi}
The condition that $v$ is almost transverse
basically means that if there is a nonconstant trajectory of $-v$ from one
critical point $p$ to a critical point $q$, then $\mbox{ind }q<\mbox{ind
}p$. Notice that for ${\rm ind}\,p<{\rm ind}\,q$ transverse intersection in 2.\ means empty
intersection and for ${\rm ind}\,p={\rm ind}\,q$ the intersection is 0-dimensional. But if
$x\in W^s(p)\cap W^u(q)$, the whole trajectory through $x$ is also in the intersection.
The existence of transverse $\omega$-gradients is given in Pajitnov \cite[Lm. 5.1]{pajiov}
which is a version of the classical Kupka-Smale theorem.\\[0.2cm]
Let $\Phi:M\times\mathbb{R}\to M$ be the flow obtained from a
vector field $v$ by integration. By a \em closed orbit \em of $v$ we mean a nonconstant map
$\gamma:S^1\to M$ with $\gamma'(x)=v(\gamma(x))$. The \em multiplicity \em
$m(\gamma)$ is the largest positive integer $m$ such that $\gamma$ factors
through an $m$-fold covering $S^1\to S^1$. Alternatively we write
$\gamma:[0,p]\to M$ with $\gamma(0)=\gamma(p)$. The number $p$ is then called
the \em period \em of $\gamma$, which we also denote by $p(\gamma)$. We say
two closed orbits are \em equivalent \em if they only differ by a rotation of
$S^1$. We denote the set of equivalence classes of closed orbits by $Cl(v)$. Notice that
$\gamma\in Cl(v)$ gives a well defined element $\{\gamma\}\in\Gamma$.\\[0.2cm]
Given $b>a\geq 0$ define $F^b_a:M\times [a,b]\to M$ by restricting $\Phi$ to
$M\times [a,b]$. The results of Section \ref{opfpt} can be used directly on
$F_a^b$. As the basepath from $v$ to $F(v,a)$ we choose $\tau(t)=\Phi(v,t)$.
It is immediate that stationary points of the flow, when viewed as fixed
points of $F_a^b$, correspond to the conjugacy class of $1_G$. For a
nontrivial closed orbit $\gamma:[0,c]\to M$ we get fixed points
$(\gamma(t),c)$ which correspond to the conjugacy class of $[\gamma]$. To see
this choose the path $\mu$ to first go from $(v,a)$ to $(v,0)$, then from
$(v,0)$ to $(\gamma(0),0)$ and then from $(\gamma(0),0)$ to
$(\gamma(0),c)$.\\[0.2cm] We have
$R(F_a^b)\in\bigoplus\limits_{\gamma\in\Gamma}H_1(C_\ast(\mathbb{Z}G,\mathbb{Z}G)_\gamma)$ and
we want to let $b$ go to infinity.
\begin{lem}Let $X$ be a compact polyhedron and $E:X\times[0,2]
\to X$ be a homotopy. Define $F,G:X\times[0,1] \to X$ by $F(x,t)=E(x,t)$ and $G(x,t)=E(x,t+1)$.\\
Then $CR(E)=CR(F)+CR(G)$. \label{l53}
\end{lem}
\begin{proof}
If $E$ is cellular, define $\Lambda:X\times I\times I\to X$ by
$\Lambda(x,t,s)=E(x,t(s+1))$ and use Proposition \ref{hprop}. In the general
case use a fine enough triangulation and a simplicial approximation for $E$
and proceed as before.
\end{proof}
From this Lemma it follows that $CR(F_0^{n+1})=CR(F_0^n)+CR(F_n^{n+1})$. To get something reasonable
for $n\to\infty$ we have to be able to somehow disregard the last term. So now assume that $\omega$ is
a Morse form and our flow comes from a weak $\omega$-gradient $v$. Let $p,q$ be critical points
of $\omega$, $x\in W^u(p)\cap W^s(q)$ and $\gamma:\mathbb{R}\to M$ the trajectory of $-v$ with
$\gamma(0)=x$. Then $\gamma$ extends to a path $\bar{\gamma}:[-\infty,\infty]\to M$ from $q$ to $p$.
\begin{defi}\em We call a loop $\delta:S^1\to M$ a \em broken closed orbit of $-v$\em, if it is a finite
concatenation of such paths $\bar{\gamma}$.
\end{defi}
\begin{rema}\label{atgiveszeta}\em If $v$ is almost transverse, no nonconstant broken closed
orbits can exist since trajectories of $-v$ between critical points lower the index.
\end{rema}
If $\gamma:[0,c]\to M$ is a nontrivial closed orbit of $-v$, we get
\[\xi(\{\gamma\})=\int\limits_\gamma\omega=\int\limits_0^c\omega(\gamma'(t))\,dt=\int\limits_0^c-\omega
(v(\gamma(t))\,dt<0\]
by Lemma \ref{smalo}. A broken closed orbit $\delta$ of $-v$ also defines a conjugacy class $\{
\delta\}\in\Gamma$ which also satisfies $\xi(\{\delta\})<0$. Then we set
\[b_\omega(-v)=\sup\,\{\,\xi(\{\delta\})\in\mathbb{R}\,|\,\delta \mbox{ is a
nonconstant broken closed orbit of }-v\}.\]
In particular the vector field $-v$ has no nonconstant
broken closed orbits if and only if $b_\omega(-v)=-\infty$. Define
\[\mathcal{O}_n=\{\gamma:[0,b]\to M\,|\,\mbox{The period }b\geq n\mbox{ and }\gamma \mbox{ is a closed orbit of }-v
\}\]
and $C_n=\sup\{c\in\mathbb{R}\,|\,-\xi([\gamma])\geq c \mbox{ for all }\gamma\in\mathcal{O}_n\}\in[0,\infty]$.
Since the $\mathcal{O}_n$ decrease as sets in $n$ we get $C_n\to C\in[0,\infty]$ as $n\to\infty$.
\begin{lem}Under the above assumptions, if
$b_\omega(-v)=-\infty$ we get $C=\infty$.\label{l54} \end{lem}
\begin{proof}
The argument is similar to Hutchings \cite[\S 3.2]{hutcth}. Assume $C<\infty$, then there exists
a sequence $\gamma_n\in\mathcal{O}_n$ with $-\xi([\gamma_n])\in[0,C]$ for all $n$.\\[0.2cm]
Choose disjoint small balls around the critical points of $\omega$ such that whenever a flowline
of $v$ leaves a ball it takes a positive time $t_0>0$ to get back into a ball. For a closed orbit
$\gamma$ let $N_\gamma$ be the number of how often the closed orbit enters (and leaves) such a ball.
Because of $t_0>0$ we get $N_\gamma<\infty$. The sequence $N_{\gamma_n}$ is bounded because otherwise
$-\xi([\gamma_n])\to\infty$. This follows because $t_0>0$ and $\omega(v)\geq\varepsilon$ outside
the balls for some $\varepsilon>0$. By passing to a subsequence we can assume that $N_{\gamma_n}$
is constant to $N$.\\[0.2cm]
Choose points $x_{n,j}$ for $j=1,\ldots\hspace{-1pt},N$ on the orbit of $\gamma_n$ away from the balls such that
there is exactly one ball on the orbit between $x_{n,j}$ and $x_{n,j+1}$ and between $x_{n,N}$
and $x_{n,1}$. Also denote by $t_{n,j}$ the time it takes from $x_{n,j}$ to $x_{n,j+1}$. We
can assume that the $x_{n,j}$ converge to $x_j\in M$ and the $t_{n,j}$ converge to $t_j\in[0,\infty]$.
Notice that $\sum_jt_{n,j}=p(\gamma_n)$.
If $t_j<\infty$ the continuity of the flow implies the existence of a flow line between $x_j$ and $x_{j+1}$
If $t_j=\infty$ there is a ``broken'' flow line from $x_j$ to $x_{j+1}$ through a critical point.
At least one of the $t_j$ has to be $\infty$ because $p(\gamma_n)\to\infty$. As a result we get
a broken closed orbit of $-v$ which contradicts $b_\omega(-v)=-\infty$.
\end{proof}
In particular for $\gamma\in\Gamma-\{\gamma(1_G)\}$ there exists an $n_\gamma>0$ such that $p_\gamma(CR(F_{n_\gamma}
^m))=0$ for all $m>n_\gamma$. The only conjugacy classes where $p_\gamma(CR(F_0^{n_\gamma}))$ can fail
to be a homology class are $\gamma(1_G)$ and conjugacy classes corresponding to fixed points of
$F(\cdot,n_\gamma)$. In particular $p_\gamma(CR(F_0^{n_\gamma}))=p_{\gamma\ast}(R(F_0^{n_\gamma}))$
is a homology class.
\begin{defi}\em Let $\omega$ be a Morse form and $v$ a weak $\omega$-gradient with
$b_\omega(-v)=-\infty$. Then we define the \em noncommutative zeta function \em
of $-v$ to be the element $\zeta(-v)\in\widehat{H\!H}_1(\mathbb{Z}G)_\xi$ which
satisfies \begin{enumerate}
\item $p_{\gamma}(\zeta(-v))=0$ for $\xi(\gamma)\geq 0$.
\item $p_{\gamma}(\zeta(-v))=p_{\gamma}(CR(F_0^{n_\gamma}))\in H_1(C_\ast(\mathbb{Z}G,\mathbb{Z}G)_\gamma)$
for $\xi(\gamma)<0$.
\end{enumerate}
\end{defi}
That $\zeta(-v)$ lies indeed in $\widehat{H\!H}_1(\mathbb{Z}G)_\xi$ follows from Lemma \ref{l53} and
\ref{l54}. By Remark \ref{atgiveszeta} $\zeta(-v)$ is defined if $v$ is almost transverse.
\begin{rema}\em There is an alternative way to decribe this zeta function. Similar to Definition \ref{dnorm}
we can define a norm $\|\cdot\|$ on $\widehat{H\!H}_1(\mathbb{Z}G)_\xi$ which turns it into a
complete metric space. Because of Lemma \ref{l54} the elements $R(F_0^n)\in H\!H_1(\mathbb{Z}G)\subset\widehat{H\!H}_1(\mathbb{Z}G)_\xi$
form a Cauchy sequence. Then $\zeta(-v)=\lim_{n\to\infty}R(F^n_0)$.
\end{rema}
The name noncommutative zeta function is motivated by the following: let $v$ be a vector field as
above that only has nondegenerate closed orbits. Here a closed orbit $\gamma$ is
\em nondegenerate \em if $\det(I-dP)\not=0$, where $P$ is a Poincar\'e map corresponding to
$\gamma$. In that case we define $\varepsilon(\gamma)\in\{1,-1\}$ to be the sign of $\det(I-dP)$.
Given a closed orbit $\gamma^m:[0,mp]\to M$ where $m$ is the multiplicity of $\gamma^m$ and hence
there is a primitive loop $\gamma:[0,p]\to M$, the conjugacy class corresponding to that orbit is
$\{\gamma^m\}$. By choosing a basepath we associate to $\gamma$ an element $[\gamma]\in G$ so that $[\gamma]^m$
represents $\{\gamma^m\}$ and then we set
\[I(\gamma^m)=[[\gamma]^m[\gamma]^{-1}\otimes[\gamma]]\in
H_1(C_\ast(\mathbb{Z}G,\mathbb{Z}G)_{\{\gamma^m\}}).\]
We do not have to worry about the basepath
because of the following
\begin{lem}\label{ezlem}
Let $g,h\in G$ and $k$ be an integer. Then
$g^k\otimes g$ is homologous to $h^{-1}g^kh\otimes h^{-1}gh$.
\end{lem}
\begin{proof}
Let $x=h^{-1}g^kh\otimes h^{-1}g\otimes h+g^kh\otimes h^{-1}\otimes g-g^{k+1}\otimes h\otimes h^{-1}$.
Then
\begin{eqnarray*}
d(x)&=&h^{-1}g^{k+1}\otimes h-h^{-1}g^kh\otimes h^{-1}gh+g^kh\otimes h^{-1}g+g^k\otimes g-g^kh\otimes h^{-1}g\\
& &+g^{k+1}h\otimes h^{-1}-g^{k+1}h\otimes h^{-1}+g^{k+1}\otimes 1-h^{-1}g^{k+1}\otimes h.
\end{eqnarray*}
Since $g^{k+1}\otimes 1$ is a boundary we get the result.
\end{proof}
In analogy with Geoghegan and Nicas \cite[\S 2B]{geonic} we form the \em Nielsen-Fuller series \em
\[\Theta(\Phi)=\sum_{\gamma\in Cl(-v)}\varepsilon(\gamma)I(\gamma)\in\widehat{H\!H}_1(\mathbb{Z}G)_\xi.\]
Here $\Phi$ denotes the flow of $-v$.
We can also define the \em eta function \em of $-v$ to be the element of $\widehat{\mathbb{Q}\Gamma}_\xi$
defined by
\[\eta(-v)(\delta)=\sum_{\gamma\in Cl(-v)\atop \{\gamma\}=\delta}\frac{\varepsilon(\gamma)}{m(\gamma)}\]
This formula for the eta function is basically taken from Pajitnov \cite{pajitn}. Commutative
eta functions already appeared in Fried \cite{fried} whose exponential is then the zeta function
of the vector field. Notice that the exponential of the noncommutative eta function is not defined
since there is no well defined multiplication in $\Gamma$ in general.\\[0.2cm]
It is easy to see that $\eta(-v)=l(\Theta(\Phi))$. Now it follows from Geoghegan and Nicas
\cite[Th.2.7]{geonic} that
\begin{equation}\label{etazeta}
\zeta(-v)=\Theta(\Phi)\mbox{ and hence }\eta(-v)=l(\zeta(-v))
\end{equation}
so $\zeta(-v)$ is a generalization of $\eta(-v)$ which is defined even when $v$ has degenerate
closed orbits. Our different convention for $R(F)$ leads to the vanishing of the `$-$' sign in
\cite[Th.2.7]{geonic}, compare Example \ref{examplos}.\\[0.2cm]
Fried \cite{fried} already defined his zeta function without the requirement of nondegenerate closed
orbits using the Fuller index \cite{fuller}. We want to show that our $\zeta(-v)$ is an appropriate
generalization in that context. Because of Lemma \ref{l54} and Fuller \cite[Th.3]{fuller} the union
$C_\gamma$ of closed orbits belonging to $\gamma\in\Gamma$ is an isolated compact set in $M\times
(0,\infty)$. Choose $n>0$ such that $C_\gamma\subset M\times[0,n-1]$ and let $C$ be the union of
closed orbits in $M\times[0,n]$. By Fuller \cite[Lm.3.1]{fuller} we can perturb the vector
field $-v$ to a vector field $-v'$ with a finite number of closed orbits in $M\times[0,n]$ and no
closed orbits in the boundary of a neighborhood of $C$. By choosing the vector field $-v'$ close
enough to $-v$ we get no further closed orbits corresponding to $\gamma$ in $M\times[0,n]$, compare
the proof of \cite[Lm.3.1]{fuller}. The straight line homotopy between $-v$ and $-v'$ induces a
homotopy between $F_0^n$ and $F'$, which gives $p_\gamma(CR(F^n_0))=p_\gamma(CR(F'))$ by Proposition
\ref{hprop}. But $l(p_\gamma(CR(F')))(\gamma)=i(C_\gamma)$, the Fuller index of $C_\gamma$, by the
remarks above. Therefore $l(p_\gamma(CR(F^n_0)))(\gamma)=i(C_\gamma)$.\\[0.2cm]
Notice that the eta function can be described as
\[\eta(-v)(\gamma)=i(C_\gamma),\]
the formula in the commutative case given by Fried \cite{fried}.
\section{The Novikov complex of a Morse form}\label{novcom}
Given a Morse form $\omega$ and a transverse $\omega$-gradient $v$
we can define the \em Novikov complex \em $C_\ast(\omega,v)$ which is in each dimension $i$ a free
$\widehat{\mathbb{Z}G}_\xi$ complex with one generator for every critical point of index $i$. Here
$\xi$ is again the homomorphism induced by $\omega$. The boundary homomorphism of $C_\ast(\omega,v)$
is based on the number of trajectories between critical points of adjacent indices. For more details
see Pajitnov \cite{pajito} or Latour \cite{latour}. This chain complex is chain homotopy equivalent
to $\widehat{\mathbb{Z}G}_\xi\otimes_{\mathbb{Z}G}C_\ast^\Delta(\tilde{M})$, where $C_\ast^\Delta
(\tilde{M})$ is the simplicial chain complex
of the universal cover $\tilde{M}$ of $M$ with respect to a smooth triangulation of $M$ lifted to
$\tilde{M}$.
The chain homotopy equivalence can be chosen so that its torsion lies in a certain subgroup of
$\overline{K}_1^G(\widehat{\mathbb{Z}G}_\xi)=K_1(\widehat{\mathbb{Z}G}_\xi)/\langle\pm [g]\,|\,g\in G\rangle$.
Given $a\in\widehat{\mathbb{Z}G}_\xi$ with $\|a\|<1$, the series $\sum_{k=0}^\infty a^n$ is a well
defined element of $\widehat{\mathbb{Z}G}_\xi$ and hence the inverse of $1-a$. Therefore
$\{1-a\in\widehat{\mathbb{Z}G}_\xi\,|\,\|a\|<1\}$ is a subgroup of $\widehat{\mathbb{Z}G}_\xi^\ast$,
the group of units of $\widehat{\mathbb{Z}G}_\xi$. We denote the image of this subgroup in
$\overline{K}_1^G(\widehat{\mathbb{Z}G}_\xi)$ by $\overline{W}$. It is proven in Pajitnov \cite{pajito}
in the rational case and in Latour \cite{latour} for the general case that there is a chain homotopy
equivalence whose torsion lies in $\overline{W}$.\\[0.2cm]
Let us specify a chain map between the completed triangulated chain complex and the Novikov
complex. We can assume that a triangulation is adjusted to $v$, i.e.\ a $k$-simplex $\sigma$ intersects
the unstable manifolds $W^u(q)$ transversely for critical points of index $\geq k$, see
\cite[\S 2.3]{schuet}. Then we define
\begin{equation}\label{formula}
\varphi(v)(\sigma)=\sum_{p\in{\rm crit}_k(\omega)}[\sigma:p]\,p
\end{equation}
where ${\rm crit}_k(\omega)$ is the set of critical points of $\omega$ having index $k$ and $[
\sigma:p]\in\widehat{\mathbb{Z}G}_\xi$ is the intersection number of a lifting of $\sigma$ to
$\tilde{M}$ with translates of the unstable manifold of a lifting of the critical point $p$.\\[0.2cm]
Let us look at the case of a rational Morse form
$\omega$ first. There is an infinite cyclic covering space $p:\bar{M}\to M$ such that $p^\ast\omega
=d\bar{f}$ is exact, namely the one corresponding to $\ker\xi$. We can also assume that $0\in\mathbb{R}$
is a regular value of $\bar{f}:\bar{M}\to\mathbb{R}$. Let $N=\bar{f}^{-1}(0)$ and $b>0$ be the number such
that $\bar{f}^{-1}(b)=tN$, where $t$ is a generator of the covering transformation group. Define
$M_N=\bar{f}^{-1}([0,b])$. Then the cobordism $(M_N,N,tN)$ is equipped with a Morse function
$\bar{f}|_{M_N}:M_N\to[0,b]$. The covering map $p$ restricted to $N$ is a diffeomorphism onto its image
and we can think of $M_N$ as a splitting of $M$ along $N$.\\[0.2cm]
For this situation Pajitnov \cite{pajirn,pajiov,pajitn} defines a condition $(\mathfrak{C}')$ for an
$\omega$-gradient $v$. For the full condition we refer the reader to these papers, but informally
it can be described as follows: the condition $(\mathfrak{C}')$ requires a Morse map $\psi$ on $N$
which gives a handle decomposition on $N$ and $tN$. The vector field $v$ which lifts to a vector
field $v'$ on $M_N$ now has to satisfy a ``cellularity condition'':
whenever $p$ is a critical point of $\bar{f}$ of index $i$, it should be the case that some
thickening of $W^s(p)$ is attached to the union of the $(i-1)$-handles of $N$ and of $M_N$.
Also a thickening of an $i$-handle in $tN$ has to flow under $-v'$ into the $i$-skeleton of $N$
and $M_N$. A symmetric condition holds for the unstable manifolds and handles of $N$.\\[0.2cm]
By Pajitnov \cite[\S 5]{pajiov}, the set of transverse $\omega$-gradients satisfying $(\mathfrak{C'})$
is $C^0$-open and dense in the set of transverse $\omega$-gradients.
Such gradients should be thought of as cellular approximations to arbitrary $\omega$-gradients.\\[0.2cm]
Now let $\rho:\tilde{M}\to M$ be the universal cover, $\tilde{f}:\tilde{M}\to\mathbb{R}$ the lifting
of $\bar{f}$ and $\tilde{N}_k=\tilde{f}^{-1}(\{k\cdot b\})$ for $k\in\mathbb{Z}$.
The handle decomposition of $N$ gives rise to sets $\tilde{V}_k^{[i]}$, $\tilde{V}_k^{(i)}$
described in \cite[\S 4.4,\S 4.5]{pajitn} such that $D_i=\bigoplus_{k\in\mathbb{Z}}
\tilde{H}_i(\tilde{V}_k^{[i]}/\tilde{V}_k^{(i-1)})$ gives a finitely generated free
$\mathbb{Z}G$ module. The topological space $\tilde{V}_k^{[i]}/\tilde{V}_k^{(i-1)}$ is in fact a
wedge of thickened $i$-spheres. The vector field $-\tilde{v}$, the lifting of $-v$ to $\tilde{M}$, induces a
map $k_i:D_i\to D_i$. We can choose a basis of $D$ by choosing lifts of handles in $\tilde{N}:=\tilde{N}_0$
and this allows us to form a matrix $A_i$ that represents $k_i$. It follows that $I-A_i$ is an
invertible matrix when viewed as a matrix over $\widehat{\mathbb{Z}G}_\xi$ and the inverse is
$\sum_{k=0}^\infty A_i^k$. It is shown in \cite[\S 4.1]{schuet} that in this situation
$\varphi(v):\widehat{\mathbb{Z}G}_\xi\otimes_{\mathbb{Z}G}
C_\ast^\Delta(\tilde{M})\to C_\ast(\omega,v)$ is a natural chain homotopy equivalence whose
torsion is given by
\[\tau(\varphi(v))=\sum_{i=0}^{n-1}(-1)^{i+1}\tau(I-A_i)\in \overline{K}_1^G(\widehat{\mathbb{Z}G}_\xi).\]
The case of irrational Morse forms is treated by approximation. The following Lemma is proven in
\cite{schuet}, see also Pajitnov \cite[\S 2B]{pajisp}.
\begin{lem}\label{approx}\cite[Lm.4.2]{schuet} For a Morse form $\omega$ and an $\omega$-gradient
$v$ there exists a rational Morse form $\omega'$ with the same set of critical points and that
agrees with $\omega$ in a neighborhood of these critical points such that $v$ is also an
$\omega'$-gradient.
\end{lem}
We denote by $\xi'$ the homomorphism induced by this rational approximation $\omega'$.
Let us compare the Novikov complexes we obtain for a Morse form $\omega$ and a rational
approximation $\omega'$ that both use the same vector field $v$. The complexes are taken over
different rings, $\widehat{\mathbb{Z}G}_\xi$ and $\widehat{\mathbb{Z}G}_{\xi'}$ respectively.
But for two critical points $p,q$ of adjacent index the elements $\tilde{\partial}(p,q)\in
\widehat{\mathbb{Z}G}_{\xi}$ and $\tilde{\partial}'(p,q)\in\widehat{\mathbb{Z}G}_{
\xi'}$ agree when viewed as elements of $\widehat{\widehat{\mathbb{Z}G}}$ since both
count the number of flowlines between $\tilde{p}$ and translates of $\tilde{q}$, and these only
depend on $v$. So we can compare chain complexes even though they are over different rings.\\[0.2cm]
We say an $\omega$-gradient $v$ satisfies the condition $(\mathfrak{AC})$, if there exists a
rational Morse form $\omega'$ such that $v$ is an $\omega'$-gradient and as such it satisfies
$(\mathfrak{C}')$. Using Lemma \ref{approx} we get $C^0$-openess and density for vector fields
satisfying $(\mathfrak{AC})$ among the transverse ones.\\[0.2cm]
Now given an $\omega$-gradient $v$ satisfying $(\mathfrak{AC})$ we can use the rational approximation
to define the matrices $A_i$ as above. It is shown in \cite[\S 4.3]{schuet} that $I-A_i$ is not
just invertible over $\widehat{\mathbb{Z}G}_{\xi'}$, but also over $\widehat{\mathbb{Z}G}_\xi$ and
the inverse is again $\sum_{k=0}^\infty A_i^k$. Furthermore \cite{schuet} shows that
$\varphi(v):\widehat{\mathbb{Z}G}_\xi\otimes_{\mathbb{Z}G}
C_\ast^\Delta(\tilde{M})\to C_\ast(\omega,v)$ is again a chain homotopy equivalence with torsion
\begin{equation}\label{torsion}
\tau(\varphi(v))=\sum_{i=0}^{n-1}(-1)^{i+1}\tau(I-A_i)\in \overline{K}_1^G(\widehat{\mathbb{Z}G}_\xi).
\end{equation}
\section{Proof of the main theorem, Part 1} \label{pofmt}
We now want to draw the connection between the torsion of the chain homotopy equivalence described
in Section \ref{novcom} and $\zeta(-v)$. Because of Pajitnov \cite{pajirn,pajitn} and our paper \cite{schuet}
we expect $\tau(\varphi(v))$ to carry the information of $\zeta(-v)$. To connect these two objects
there is the \em Dennis trace \em homomorphism $DT:K_1(R)\to H\!H_1(R)$ defined by
$DT(\alpha)=[\mbox{trace }A^{-1}\otimes A]$, where $\alpha\in K_1(R)$ is represented by the matrix
$A$. Here $R$ is a ring with unit. It is elementary that $DT$ is a well defined homomorphism. For
more on the Dennis trace see Igusa \cite[\S 1]{igusa}.\\[0.2cm]
Now we choose $R=\widehat{\mathbb{Z}G}_\xi$ and the composition $\theta\circ DT$ is a homomorphism
\\$K_1(\widehat{\mathbb{Z}G}_\xi)\to \widehat{H\!H}_1(\mathbb{Z}G)_\xi$. Our torsion $\tau(\varphi(v))$
is an element of $\overline{W}\subset\overline{K}_1^G(\widehat{\mathbb{Z}G}_\xi)$. We can also
define a subgroup $W\subset K_1(\widehat{\mathbb{Z}G}_\xi)$ as the image of $\{1-a\,|\,\|a\|<1\}$
in $K_1(\widehat{\mathbb{Z}G}_\xi)$.
\begin{lem}
The projection $p: K_1(\widehat{\mathbb{Z}G}_\xi)\to\overline{K}_1^G(\widehat{\mathbb{Z}G}_\xi)$
restricted to $W$ induces an isomorphism $W\to\overline{W}$.\footnote{In the case of a rational
Novikov ring this also follows from Pajitnov and Ranicki \cite[Cor.0.1]{pajran}}
\end{lem}
\begin{proof}
Look at the composition
\[ K_1(\widehat{\mathbb{Z}G}_\xi)\stackrel{DT}{\longrightarrow}H\!H_1(\widehat{\mathbb{Z}G}_\xi)
\stackrel{\theta}{\longrightarrow}\widehat{H\!H}_\ast(\mathbb{Z}G)_\xi\stackrel{p_{\gamma(1)}}
{\longrightarrow}H_1(C_\ast(\mathbb{Z}G,\mathbb{Z}G)_{\gamma(1)})=H_1(G)=G_{ab}.\]
It is easy to see that the image of $\tau(1-a)\in W$ under this composition in $G_{ab}$ is 0. Denote
the image of $\pm G$ in $K_1(\widehat{\mathbb{Z}G}_\xi)$ by $\bar{G}$. Then $\tau(\pm g)\in\bar{G}$
gets mapped to $g[G,G]$, compare Geoghegan and Nicas \cite[\S 6A]{geonia}. So $W\cap\bar{G}\subset
\{\tau(\pm 1)\}$. To see that $\tau(-1)\notin W$ look at the composition
$K_1(\widehat{\mathbb{Z}G}_\xi)\stackrel{\varepsilon_\ast}{\longrightarrow}
K_1(\widehat{\mathbb{Z}H_1(G)}_{\bar{\xi}})\stackrel{\det}{\longrightarrow}
(\widehat{\mathbb{Z}H_1(G)}_{\bar{\xi}})^\ast$.
\end{proof}
Therefore we get a homomorphism $\mathfrak{DT}=
\theta\circ DT\circ(p|_W)^{-1}:\overline{W}\to\widehat{H\!H}_1(\mathbb{Z}G)_\xi$ and we want to
compare $\mathfrak{DT}(\tau(\varphi(v)))$ with $\zeta(-v)$.\\[0.2cm]
So given a Morse form $\omega$ we denote by $\mathcal{GA}(\omega)$ the set of $\omega$-gradients
satisfying $(\mathfrak{AC})$. The theorem we can prove now reads
\begin{theo}\label{mtheo}
Let $\omega$ be a Morse form on a smooth connected closed manifold $M^n$. Let $\xi:G\to\mathbb{R}$
be induced by $\omega$ and let $C^\Delta_\ast(\tilde{M})$ be the simplicial $\mathbb{Z}G$ complex
coming from a smooth triangulation of $M$. For every $v\in\mathcal{GA}(\omega)$ there is a natural
chain homotopy equivalence $\varphi(v):\widehat{\mathbb{Z}G}_\xi\otimes_{\mathbb{Z}G}C_\ast^\Delta
(\tilde{M})\to C_\ast(\omega,v)$ given by (\ref{formula}) whose torsion $\tau(\varphi(v))$ lies
in $\overline{W}$ and satisfies
\[\mathfrak{DT}(\tau(\varphi(v)))=\zeta(-v).\]
\end{theo}
\begin{rema}\em
This theorem is a generalization of Pajitnov \cite[Main Th.]{pajitn} and of our paper \cite[Th.4.5]{schuet} in
that $\zeta(-v)$ is a generalization of $\eta(-v)$ and the condition that the vector fields $v$
only have nondegenerate closed orbits is dropped.
\end{rema}
To prove Theorem \ref{mtheo} we have to show for every $\gamma\in\Gamma$ that $p_{\gamma\ast}
\circ\mathfrak{DT}(\tau(\varphi(v)))=p_{\gamma\ast}(\zeta(-v))$.
For this we have to compare the Hochschild chains that represent both sides of the
equation and show that they are homologous. After we bring the chains into a certain form the last
comparison will follow from the fact that the Lefschetz number can be computed by the fixed point
index. This idea is basically given in Pajitnov \cite[\S 8]{pajirn} but since
$\widehat{H\!H}_1(\mathbb{Z}G)_\xi$ is more delicate than
$\widehat{\mathbb{Q}\Gamma}_\xi$ the proof will be slightly more
involved.\\[0.2cm]
Because of (\ref{torsion}) let us look at $\mathfrak{DT}(\tau(I-A))$, where
$A$ is an $l\times l$ matrix over $\mathbb{Z}G$ such that $I-A$ is invertible
over $\widehat{\mathbb{Z}G}_\xi$ with inverse $I+A+A^2+\ldots$. We denote the
entries of $A$ by $A_{ij}\in\mathbb{Z}G$. We can consider $\tau(I-A)\in
K_1(\widehat{\mathbb{Z}G}_\xi)$. Then
\[DT(\tau(I-A))=[\mbox{trace}(\sum_{k=1}^\infty A^{k-1}\otimes
(I-A))]=-[\mbox{trace}( \sum_{k=1}^\infty A^{k-1}\otimes A)]\in
H\!H_1(\widehat{\mathbb{Z}G}_\xi).\] Passing to
$\widehat{H\!H}_1(\mathbb{Z}G)_\xi$ by $\theta$ allows us to move the series
out of the tensor, so let us look at $\mbox{trace}(A^{k-1}\otimes A)\in
C_1(\mathbb{Z}G,\mathbb{Z}G)$. We have \begin{equation}\label{traceA}
\mbox{trace}(A^{k-1}\otimes A)=\sum_{i_1,\ldots\hspace{-1pt},i_k=1}^l
A_{i_1i_2}A_{i_2i_3}\cdots A_{i_{k-1}i_k} \otimes A_{i_ki_1}. \end{equation}
We need to bring the 1-chain (\ref{traceA}) into a different form.
\begin{lem}\label{homtrace} Let $k,r,i_1,\ldots\hspace{-1pt},i_k$ be positive
integers and $A_{i_s,i_t}\in\mathbb{Z}G$ for
$s,t\in\{1,\ldots\hspace{-1pt},k\}$. Denote $A_\ast=A_{i_1i_2}\cdots
A_{i_ki_1}$, then $A_\ast^{r-1}A_{i_1i_2} \cdots A_{i_{k-1}i_k}\otimes
A_{i_ki_1}+ A_{i_2i_3}\cdots A_{i_ki_1}A_\ast^{r-1}\otimes A_{i_1i_2}
+\ldots+A_{i_ki_1}A_\ast^{r-1}A_{i_1i_2}\cdots A_{i_{k-2}i_{k-1}}\otimes
A_{i_{k-1}i_k}$ is homologous to $A_\ast^{r-1}\otimes A_\ast$. \end{lem}
\begin{proof} Look at the 2-chain $A_\ast^{r-1}\otimes A_{i_1i_2}\otimes
A_{i_2i_3}\cdots A_{i_ki_1}+A_\ast^{r-1} A_{i_1i_2}\otimes A_{i_2i_3}\otimes
A_{i_3i_4}\cdots A_{i_ki_1}+\ldots +A_\ast^{r-1}A_{i_1i_2}\cdots
A_{i_{k-2}i_{k-1}}\otimes A_{i_{k-1}i_k}\otimes A_{i_ki_1}$. It is
straightforward to check that its boundary gives the result. \end{proof}
\begin{proof}[Proof of Theorem \ref{mtheo}]
Since $v\in\mathcal{GA}(\omega)$, there exists a rational Morse form $\omega'$ such that $v$
satisfies condition $(\mathfrak{C'})$ with respect to $\omega'$. Let $\xi'$ be the homomorphism
induced by $\omega'$, $\bar{M}$ the infinite cyclic covering space corresponding to $\ker\xi'$
and $\bar{f}:\bar{M}\to\mathbb{R}$ a smooth function with $0\in\mathbb{R}$ a regular value and $
p^\ast\omega'=d\bar{f}$. Furthermore choose $b>0$ such that $M_N=\bar{f}^{-1}([0,b])$ is a
splitting of $M$ along $N$. We can assume that $v$ satisfies condition $(\mathfrak{C}')$ with
respect to this splitting. Hence there is a Morse function $\psi:N\to[0,1]$ and on $N_k=\bar{f}
^{-1}(\{kb\})$ for all $k\in\mathbb{Z}$ which is ordered in the sense of Milnor \cite[Def.4.9]{milnhc}.
Then we get filtrations of $N$ and $N_k$ by
\[V^{(i)}_k=\psi^{-1}([0,\alpha_{i+1}])\subset N_k,\]
where $\alpha_{i+1}$ is a real number bigger than the image of critical points of index $i$ under
$\psi$ and smaller than the image of critical points of index $i+1$ under $\psi$, and
\[V^{[i]}_k=V^{(i-1)}_k\cup\mbox{ thickenings of the stable manifolds of
critical points of index }i \subset V^{(i)}_k.\]
If $\gamma$ is a closed orbit of $-v$, it lifts to a trajectory $\bar{\gamma}:\mathbb{R}\to\bar{M}$
such that $\bar{f}\circ\bar{\gamma}:\mathbb{R}\to\mathbb{R}$ is bijective. Assume $\bar{\gamma}(0)\in
N$ and let $p>0$ be the prime period of $\gamma$, i.e.\ the smallest $p>0$ such that $\gamma(p)=
\gamma(0)$. There is a negative integer $k$ such that
$\bar{\gamma}(p)\in N_k$. Furthermore $t^k\bar{\gamma}(0)=\bar{\gamma}(p)$, where $t$ is the
generator of the covering transformation group that satisfies $tN_0=N_1$. Let $i$ be a number
such that $\bar{\gamma}(0)\in V^{(i)}_0-V^{(i-1)}_0$. It follows from the definition of condition
$(\mathfrak{C}')$ that $\bar{\gamma}(\mathbb{R})\cap N_k\subset V^{[i]}_k$ for all $k\in\mathbb{Z}$,
see Pajitnov \cite{pajirn,pajitn}. In fact the points are in the interior of the handle.
Therefore we get a partition of $Cl(-v)$ into sets $Cl(-v;i)$ consisting of closed orbits passing
through $V_k^{(i)}-V^{(i-1)}_k$.\\[0.2cm]
We get the chain homotopy equivalence $\varphi(v):\widehat{\mathbb{Z}G}_\xi\otimes_{\mathbb{Z}G}
C_\ast^\Delta(\tilde{M})\to C_\ast(\omega,v)$ with torsion (\ref{torsion}) from Section \ref{novcom}.
We show that the $i$-th summand $(-1)^{i+1}\tau(I-A_i)$ carries the information of $Cl(-v;i)$, so
let us fix $i$ and for ease of notation denote the matrix $A_i$ by $A$. To describe the entries
$A_{jk}$ of $A$ we need the universal cover $p:\tilde{M}\to\bar{M}$. For $X\subset\bar{M}$ let
$\tilde{X}=p^{-1}(X)$. The map $\bar{f}$ induces a map $\tilde{f}:\tilde{M}\to\mathbb{R}$. The matrix
$A$ described in Section \ref{novcom} comes in fact from a map $(-\tilde{v})^\to:
\tilde{V}_0^{[i]}/\tilde{V}^{(i-1)}_0\to\tilde{V}_{-1}^{[i]}/\tilde{V}_{-1}^{(i-1)}$
induced by $-\tilde{v}$. Now $\tilde{V}_0^{[i]}/\tilde{V}^{(i-1)}_0$ is a wedge of thickened $i$-spheres
and $\coprod_{k\in\mathbb{Z}}\tilde{V}_k^{[i]}/\tilde{V}^{(i-1)}_k$ has a natural $G$-action and
modulo $G$ it consists of as many spheres as $\psi$ has critical points of index $i$. To get the
matrix $A$ we have to lift the handles of the critical points to $\tilde{N}=\tilde{N}_0$. Hence
we pick spheres in $\tilde{V}_0^{[i]}/\tilde{V}^{(i-1)}_0$ that we denote by $\sigma_j$. So if we
denote the composition
\[\sigma_j\hookrightarrow \tilde{V}_0^{[i]}/\tilde{V}^{(i-1)}_0\stackrel{(-\tilde{v})^\to}
{\longrightarrow}\tilde{V}_{-1}^{[i]}/\tilde{V}^{(i-1)}_{-1}\stackrel{r}{\longrightarrow}g\sigma_k\]
by $\sigma^g_{jk}$, where the last map just retracts every sphere other than $g\sigma_k$ to
the wedge point, then $A_{jk}(g)$ is the degree of $\sigma^g_{jk}$.\\[0.2cm]
Now we look at $\mbox{trace}(A^{k-1}\otimes A)$. This term contains information about the map
\\$((-\tilde{v})^\to)^k:\tilde{V}_0^{[i]}/\tilde{V}^{(i-1)}_0\to\tilde{V}_{-k}
^{[i]}/\tilde{V}^{(i-1)}_{-k}$.\\[0.2cm]
We say $(i_1,\ldots\hspace{-1pt},i_k),(j_1,\ldots\hspace{-1pt},j_k)\in\{1,\ldots\hspace{-1pt},l\}^k$ are
\em equivalent \em if they differ only by a rotation and denote by $S$ the set of equivalence classes. Then
\begin{eqnarray*}
\mbox{trace}(A^{k-1}\otimes A)&=&\sum_{i_1,\ldots\hspace{-1pt},i_k=1}^l A_{i_1i_2}\cdots A_{i_{k-1}i_k}
\otimes A_{i_ki_1}\\
&=&\sum_{[x]\in S}\sum_{(i_1,\ldots\hspace{-1pt},i_k)\in[x]}A_{i_1i_2}\cdots A_{i_{k-1}i_k}\otimes A_{i_ki_1}
\end{eqnarray*}
Fix $[x]\in S$. The order of $[x]$ divides $k$. Let $q$ be the order of $[x]$ and let $r$ be so that
$qr=k$. If $(i_1,\ldots\hspace{-1pt},i_k)\in[x]$, then $(i_1,\ldots\hspace{-1pt},i_k)=(i_1,\ldots\hspace{-1pt},i_q,\ldots\hspace{-1pt},i_1,\ldots\hspace{-1pt},i_q)$.
\\[0.2cm]
Let us first look at $r=1$. By Lemma \ref{homtrace} $\sum\limits_{(i_1,\ldots\hspace{-1pt},i_k)\in[x]}
A_{i_1i_2}\cdots A_{i_{k-1}i_k}\otimes A_{i_ki_1}$ is homologous to $1\otimes A_{i_1i_2}\cdots
A_{i_ki_1}$. For $g_j\in\mbox{supp }A_{i_ji_{j+1}}$, $A_{i_1i_2}(g_1)\cdots A_{i_ki_1}(g_k)$ is the degree of
the map
\[(g_1\cdots g_{k-1}\sigma^{g_k}_{i_ki_1})\circ\ldots\circ(g_1\sigma^{g_2}_{i_2i_3})\circ\sigma^{g_1}
_{i_1i_2}:\sigma_{i_1}\to g_1\cdots g_k\sigma_{i_1}.\]
Let $\chi:\sigma_{i_1}\to\sigma_{i_1}$ be the composition of that map with $(g_1\cdots g_k)^{-1}:g_1\cdots g_k\sigma_{i_1}\to
\sigma_{i_1}$. Then a fixed point other than the basepoint of $\chi$ corresponds to a closed orbit
of $-v$. Notice that a closed orbit arising this way has multiplicity 1, since it passes through
the spheres $\sigma_j$ without a repeating pattern (recall $r=1$). Furthermore $(i_1,\ldots\hspace{-1pt},i_k)$
describes the cells through which the closed orbit passes.\\[0.2cm]
To continue we now look at the special case where the vector field $-v$ only has nondegenerate
closed orbits. Then the fixed points of $\chi$ are isolated. Notice that the basepoint is a fixed
point and its index is 1, since $\chi$ is constant near the basepoint. Now the Lefschetz number of
$\chi$ satisfies
\begin{equation}\label{lefsch}
L(\chi)=1+(-1)^iA_{i_1i_2}(g_1)\cdots A_{i_ki_1}(g_k)=\sum\mbox{ fixed point indices}
\end{equation}
so if $Cl(-v;i;i_1,\ldots\hspace{-1pt},i_k;g_1,\ldots\hspace{-1pt},g_k)$ is the subset of $Cl(-v;i)$ consisting of closed
orbits following the pattern $\sigma_{i_1}\to g_1\sigma_{i_2}\to\ldots\to g_1\cdots g_k\sigma_{i_1}$,
then
\[(-1)^{i}A_{i_1i_2}(g_1)\cdots A_{i_ki_1}(g_k)=\sum_{\gamma\in Cl(-v;i;i_1,\ldots\hspace{-1pt},i_k;g_1,\ldots
,g_k)}\varepsilon(\gamma).\]
Recall that $\varepsilon(\gamma)$ is the fixed point index of the Poincar\'e map and hence the
fixed point index of the corresponding fixed point of $\chi$. Notice that $\gamma\in Cl(-v;i;i_1,
\ldots,i_k;g_1,\ldots\hspace{-1pt},g_k)$ contributes the 1-chain $\varepsilon(\gamma) 1\otimes g_1\cdots g_k$
to $\Theta(\Phi)=\zeta(-v)$ and $g_1\cdots g_k=[\gamma]$.
By looking at all combinations $g_j\in\mbox{supp }A_{i_ji_{j+1}}$ we get that
\[(-1)^i\sum_{(i_1,\ldots\hspace{-1pt},i_k)\in[x]}A_{i_1i_2}\cdots A_{i_{k-1}i_k}\otimes A_{i_ki_1}\sim\sum_{
\gamma\in Cl(-v;i;[x])}\varepsilon(\gamma)1\otimes [\gamma]\]
where $\sim$ means homologous and $Cl(-v;i;[x])$ is the subset of $Cl(-v;i)$ consisting of closed
orbits following a pattern $\sigma_{i_1}\to g_1\sigma_{i_2}\to\ldots\to g_1\cdots g_k\sigma_{i_1}$
for $(i_1,\ldots\hspace{-1pt},i_k)\in[x]$ and some $g_j\in\mbox{supp }A_{i_ji_{j+1}}$.\\[0.2cm]
Now look at $r>1$. Then $A_{i_1i_2}\cdots A_{i_{k-1}i_k}\otimes A_{i_ki_1}=(A_{i_1i_2}\cdots A_{i_q
i_1})^{r-1}A_{i_1i_2}\cdots A_{i_{q-1}i_q}\otimes A_{i_qi_1}$ and
\[\sum_{(i_1,\ldots\hspace{-1pt},i_k)\in[x]}A_{i_1i_2}\cdots A_{i_{k-1}i_k}\otimes A_{i_ki_1}\sim
(A_{i_1i_2}\cdots A_{i_qi_1})^{r-1}\otimes A_{i_1i_2}\cdots A_{i_qi_1}\]
by Lemma \ref{homtrace}. Let $g_j\in\mbox{supp }A_{i_ji_{j+1}}$ for $j=1,\ldots\hspace{-1pt},k$.
 Note that $i_j=i_{j+q}$ but we can have $g_j\not= g_{j+q}$. We also have $g_j\in\mbox{supp }A_{i_j
i_{j+1}}$ implies $\xi'(g_j)=-b$, so we restrict our attention to $G_{-1}=(\xi')^{-1}(\{-b\})$.
Again $A_{i_1i_2}(g_1)\cdots A_{i_ki_1}(g_k)$ is the degree of $\chi$ which is defined just as in
the case $r=1$, but this time fixed points can correspond to closed orbits $\gamma$ with multiplicity
$m(\gamma)>1$.\\[0.2cm]
We say $(h_1,\ldots\hspace{-1pt},h_k),(h'_1,\ldots\hspace{-1pt},h'_k)\in(G_{-1})^k$ are equivalent,
if they differ only by a rotation of $jq$ elements, where $j$ is an integer, e.g. for $x,y$ different
$(x,y,x,y,x,y)$ and $(y,x,y,x,y,x)$ are equivalent for $q=3$, but not for $q=2$. Denote the set
of equivalence classes by $T$.\\[0.2cm]
We now have
\begin{eqnarray*}
(A_{i_1i_2}\cdots A_{i_qi_1})^{r-1}\otimes A_{i_1i_2}\cdots A_{i_qi_1}&=&
\sum_{g_1,\ldots\hspace{-1pt},g_k\in G_{-1}}n_{g_1}\cdots n_{g_k} g_1\cdots g_{k-q}\otimes g_{k-q+1}
\cdots g_k\\
&=&\sum_{[y]\in T}\sum_{(g_1,\ldots\hspace{-1pt},g_k)\in[y]} n_{g_1}\cdots n_{g_k} g_1\cdots g_{k-q}
\otimes g_{k-q+1}\cdots g_k.
\end{eqnarray*}
Note that $n_{g_1}\cdots n_{g_k}$ only depends on $[y]$, so we denote this by $n_{[y]}$. Fix $[y]$
and let $s$ be the order of $[y]$. Then $s$ divides $r$, so let $s\cdot p=r$ and we get for
$(g_1,\ldots\hspace{-1pt},g_k)\in[y]$ that $(g_1,\ldots\hspace{-1pt},g_k)=
(g_1,\ldots\hspace{-1pt},g_{qp},\ldots\hspace{-1pt},g_1,\ldots\hspace{-1pt},g_{qp})$, where
$(g_1,\ldots\hspace{-1pt},g_{qp})$ repeats $s$ times. By Lemma \ref{homtrace} we get
\[\sum_{(g_1,\ldots\hspace{-1pt},g_k)\in[y]} g_1\cdots g_{k-q}\otimes g_{k-q+1}\cdots g_k\sim
(g_1\cdots g_{qp})^{s-1}\otimes g_1\cdots g_{qp}.\]
As mentioned above, fixed points of $\chi$ correspond to closed orbits of $-v$. Let $\gamma$ be a
closed orbit coming from a fixed point of $\chi$ which is not the basepoint. The multiplicity of
$\gamma$ divides $s$, let $m'(\gamma)$ be the number with $m(\gamma)m'(\gamma)=s$. The closed orbit
$\gamma$ then provides $m'(\gamma)$ different fixed points of $\chi$, all with the same fixed point
index. Furthermore $I(\gamma)$ is represented by $(g_1\cdots g_{qp})^{s-m'(\gamma)}\otimes(g_1
\cdots g_{qp})^{m'(\gamma)}$ which is homologous to $m'(\gamma)\cdot(g_1\cdots g_{qp})^{s-1}
\otimes g_1\cdots g_{qp}$ by Lemma \ref{homtrace}. Denote a representing chain of $I(\gamma)$ by $I'(\gamma)$. Now we can
basically proceed as in the case $r=1$,
\[(-1)^in_{[y]}=\sum_{\gamma\in Cl(-v;i;[y])}\varepsilon(\gamma) I'(\gamma)\]
where $Cl(-v;i;[y])$ consists of closed orbits following a pattern $\sigma_{i_1}\to g_1\sigma_{i_2}
\to\ldots\to g_1\cdots g_k\sigma_{i_1}$ for $(g_1,\ldots\hspace{-1pt},g_k)\in[y]$ and the right
side is the contribution of those closed orbits to $\zeta(-v)$. Summing over all $[y]\in T$ gives
\[(-1)^i(A_{i_1i_2}\cdots A_{i_qi_1})^{r-1}\otimes A_{i_1i_2}\cdots A_{i_qi_1}\sim\sum_{\gamma
\in Cl(-v;i;[x])}\varepsilon(\gamma)I'(\gamma).\]
Therefore we get that
\[-(-1)^{i+1}\mbox{trace}(A^{k-1}\otimes A)\sim\sum_{\gamma\in Cl(-v;i;k)}\varepsilon(\gamma)I'
(\gamma).\]
Here $Cl(-v;i;k)$ consists of those $\gamma:[0,p(\gamma)]\to M\in Cl(-v;i)$ that satisfy $\bar{\gamma}
(0)\in N_0$ and $\bar{\gamma}(p(\gamma))\in N_{-k}$ for a lift $\bar{\gamma}$ of $\gamma$ to
$\bar{M}$.\\[0.2cm]
Because of (\ref{etazeta}) and (\ref{torsion})
summing over all $k$ and all $i$ gives the desired result for vector fields that only have nondegenerate
closed orbits.\\[0.2cm]
Now we have to allow degenerate closed orbits for our vector field $v$. A key fact in proving the
theorem for vector fields with only nondegenerate closed orbits was the equality (\ref{lefsch})
which still holds in the general case. The fixed point index has to be taken in a more general
sense, see Brown \cite[Ch.4]{brown}. What needs to be shown is that the fixed point index contains
the right information for $\zeta(-v)$.\\[0.2cm]
So fix $k$ and look at $\mbox{trace}(A^{k-1}\otimes A)$. The matrix $A$ comes from the flow of
$-v$. Recall the maps $\sigma^g_{jm}$ which have $A_{jm}(g)$ as degree. By abuse of notation we
denote the $i$-handle in $V_0^{[i]}\subset N_0$ corresponding to the thickened $i$-sphere $\sigma_j$
also by $\sigma_j$. If $x\in\sigma_j$ is a point whose $-v$-trajectory leads into a critical point
of $\bar{f}$ before it reaches $N_{-1}$, then $\sigma^g_{jm}(x)=\ast$, the basepoint of $g\sigma_m
\subset \tilde{V}^{[i]}_{-1}/\tilde{V}^{(i-1)}_{-1}$. By the definition of $(-\tilde{v})
^\to$ we get that points near $x$ will also be mapped to $\ast$ under $\sigma^g_{jm}$.
Also points on a closed orbit have to be in the interior of the handle by the definition of
condition $(\mathfrak{C}')$. Points in $\sigma_j$ that do not get mapped to $\ast$ under $\sigma^g_{jm}$
have trajectories avoiding critical points of $\bar{f}$ between $N_{-1}$ and $N_0$. A compactness
argument gives that the information to define $\sigma^g_{jm}$ and up to $k$ compositions of these
maps is contained in a finite piece of the flow. Let $F:M\times [0,s]\to M$ be such a piece.
The fixed point set of $F$ does not consist of finitely many circles in general (aside the fixed
points of $F_0$, $F_s$ and stationary points). But by transversality we can change $F$ to a homotopy
$F'$ whose fixed points are finitely many circles. By changing only near the fixed points of $F$
we can use the homotopy between $F$ and $F'$ to get a map $(\sigma^g_{jm})'$ homotopic to
$\sigma^g_{jm}$ which has finitely many fixed points, each corresponding to a circle of $\mbox{Fix}
(F')$.\\[0.2cm]
Now $R(F)=R(F')$ and by Geoghegan and Nicas \cite[Th.1.10]{geonia} $R(F')$ can be computed by
the transverse intersection invariant $\Theta(F')$, see also \cite[Th.1.15]{geonic} for an easy
interpretation of $\Theta(F')$ in terms of 1-chains. But $\Theta(F')$ gives us the right
comparison with $\tau(\varphi(v))$ just as the Nielsen Fuller series did in the case of nondegenerate
orbits only. This finishes the proof of Theorem \ref{mtheo}.
\end{proof}
\section{Properties of the zeta function}
We want to remove the cellularity condition in Theorem \ref{mtheo}. To do this we will show
that the zeta function and the torsion depend continuously on the vector field. The statement
will then follow from the density of $\mathcal{GA}(\omega)$. This section is about the
zeta function and the remaining sections will deal with the torsion.
\begin{defi}\em Let $(v_t)_{t\in[0,1]}$ be a smoothly varying one parameter family of weak
$\omega$-gradients and $R\in[-\infty,0)$. We say $(v_t)$ is \em $R$-controlled\em, if
$b_\omega(-v_t)\leq R$ for all $t\in[0,1]$. \end{defi}
\begin{prop}\label{contzeta}
Let $R\in[-\infty,0)$ and $(v_t)_{t\in[0,1]}$ be an $R$-controlled one parameter family of weak
$\omega$-gradients such that $b_\omega(-v_0)=b_\omega(-v_1)=-\infty$. Then
$p_\gamma(\zeta(-v_0))=p_\gamma( \zeta(-v_1))$ for every $\gamma\in\Gamma$
with $\xi(\gamma)>R$. \end{prop}
\begin{proof}
We need an argument similar to the proof of Lemma \ref{l54}. Define
\[\mathcal{O}_n=\{c:[0,b]\to M\,|\, b\geq n \mbox{ and }c\mbox{ is a closed orbit of }-v_t
\mbox{ for some }t\in[0,1]\}\]
and $C_n=\sup\,\{x\in\mathbb{R}\,|\,-\xi([c])\geq x\mbox{ for all }c\in\mathcal{O}_n\}\in[0,\infty]$.
Again we get $C_n\to C\in[0,\infty]$ as $n\to\infty$. We claim that $C\geq -R$.\\[0.2cm]
So let $\varepsilon>0$ and assume $C\leq -(R+\varepsilon)$. Now we proceed as in the proof of Lemma
\ref{l54} to get a broken closed orbit $\delta$ of $-v_t$ for some $t\in[0,1]$. But by continuity
we get $-\xi(\{\delta\})\leq-(R+\varepsilon)$ contradicting the fact that the one parameter family
is $R$-controlled.\\[0.2cm]
Now let $\gamma\in\Gamma$ satisfy $\xi(\gamma)>R$ and let $F_t$ be the flow of $-v_t$. Then
$F:M\times\mathbb{R}\times [0,1]\to M$ is a smooth homotopy of flows and by the argument above
there exists an $n>0$ such that $F|_{M\times[0,n]\times[0,1]}$ contains all closed orbits $c$ with
$\{c\}=\gamma$. By Proposition \ref{hprop} we get
\[p_\gamma(\zeta(-v_0))-p_\gamma(\zeta(-v_1))=-p_\gamma(CR(U)),\]
where $U:M\times[0,1]\to M$ is given by $U(x,t)=F(x,n+1,t)$. But $p_\gamma(CR(U))=0$, since $U$ has
no fixed points corresponding to $\gamma$ as such fixed points would give a closed orbit of
some $-v_t$ with period $n+1$ corresponding to $\gamma$.
\end{proof}
\begin{coro}
Let $\omega$ be a Morse form and $v_0,v_1$ be weak $\omega$-gradients such that there exists a
$(-\infty)$-controlled one parameter family joining them. Then $\zeta(-v_0)=\zeta(-v_1)$.
\end{coro}
If $D_\ast$ is a free finitely generated acyclic complex over a ring $R$, denote by $\tau(D_\ast)$
its torsion.
\begin{coro}\label{constnc}
Let $\omega$ be a closed 1-form without critical points and $v$ an $\omega$-gradient. Then
\[\zeta(-v)=-\mathfrak{DT}(\tau(\widehat{\mathbb{Z}G}_\xi\otimes_{\mathbb{Z}G}
C_\ast^\Delta(\tilde{M}))).\]
In particular $\zeta(-v)$ does not depend on $v$.
\end{coro}
\begin{proof}
We have $b_\omega(-v)=-\infty$ for any $\omega$-gradient, since there are no
critical points. Let $w\in \mathcal{GA}(\omega)$. By Theorem \ref{mtheo} we
have $\zeta(-w)=\mathfrak{DT}(\tau(\varphi))$, where
$\varphi:\widehat{\mathbb{Z}G}_\xi\otimes_{\mathbb{Z}G}C_\ast^\Delta(\tilde{M})\to C_\ast(\omega,v)=0$ is the zero map. But for any chain homotopy equivalence $\psi:D_\ast\to E_\ast$
between acyclic complexes we have $\tau(\psi)=\tau(E_\ast)-\tau(D_\ast)$. Hence
$\tau(\varphi)=-\tau(\widehat{\mathbb{Z}G}_\xi\otimes_{\mathbb{Z}G}C_\ast^\Delta(\tilde{M}))$.
Now $\zeta(-v)=\zeta(-w)$ since the one parameter family $v_t=tv+(1-t)w$ is $(-\infty)$-controlled.
\end{proof}
We now want to show that Proposition \ref{contzeta} is still useful when we allow critical points.\\[0.2cm]
Let $(W;N_0,N_1)$ be a cobordism and $f:W\to[a,b]$ a Morse function on $W$ in the sense of Milnor
\cite[Def.2.3]{milnhc}. Just as for closed 1-forms we can define $f$-gradients and the notion
of transverse and almost transverse $f$-gradients. We say a weak $f$-gradient is almost transverse
if for two different critical points $p,\,q$ of $f$ with ${\rm ind}\,p\leq{\rm ind}\,q$ we have
$W^s(p)\cap W^u(q)=\emptyset$.
\begin{lem}\label{pajlem}
Let $v$ be an almost transverse $f$-gradient on $W$. Then any weak $f$-gradient $w$ sufficiently
close to $v$ in the $C^0$-topology is also almost transverse.
\end{lem}
\begin{proof}
Since $v$ is an almost transverse $f$-gradient, we can rearrange $f$ to a Morse function $\phi$
which is self-indexed in the sense of Milnor \cite[Df.4.9]{milnhc}, near the critical points $\phi-f$
is constant and such that $v$ is also a
$\phi$-gradient. By Pajitnov \cite[Lm.2.74]{pajiov} every weak $f$-gradient $w$ close enough to
$v$ is also a weak $\phi$-gradient. But such a vector field is almost transverse, since for
critical points $p\not=q$ with ${\rm ind}\,p\leq{\rm ind}\,q$ we have $\phi(p)\leq\phi(q)$ and so
$W^s(p,w)\cap W^u(q,w)=\emptyset$.
\end{proof}
Let $\mathcal{G}_{at}(\omega)$ be the set of almost transverse $\omega$-gradients together with
the $C^0$-topology.
\begin{theo}\label{zetacont}
Let $\omega$ be a Morse form on the connected closed smooth manifold $M$. Then $\zeta:\mathcal{G}_{at}
(\omega)\to\widehat{H\!H}_1(\mathbb{Z}G)_\xi$ sending $v$ to $\zeta(-v)$ is continuous.
\end{theo}
\begin{proof}
Let $v\in\mathcal{G}_{at}(\omega)$. Given $R<0$ we need to find a neighborhood $U(v)$ of $v$ in
$\mathcal{G}_{at}(\omega)$ such that for all $\gamma\in\Gamma$ with $\xi(\gamma)\geq R$ we have
$p_\gamma(\zeta(-v))=p_\gamma(\zeta(-w))$ for all $w\in\mathcal{G}_{at}(\omega)$.\\[0.2cm]
Assume first that $\omega$ is rational. As in Section \ref{novcom} we have the infinite cyclic
covering space $p:\bar{M}\to M$ and a smooth function $\bar{f}:\bar{M}\to\mathbb{R}$ with
$p^\ast\omega=d\bar{f}$ and $0\in\mathbb{R}$ as a regular value. Let $b>0$ as in Section \ref{novcom}.
For any two regular values $a_1<a_2$ $\bar{f}|_{\bar{f}^{-1}[a_1,a_2]}$ is a Morse function on
the cobordism $W_{a_1,a_2}$. The vector field $v$ lifts to a transverse $\bar{f}$-gradient $\bar{v}$
on $\bar{M}$ whose restriction to $W_{a_1,a_2}$ is an $\bar{f}|$-gradient. If we choose another
$\omega$-gradient $w$ close to $v$, then its lift $\bar{w}$ will also be close to $\bar{v}$.
Since $\bar{v}$ is almost transverse on $W_{a_1,a_2}$, an $\omega$-gradient $w$ close enough to
$v$ will lift to an $\bar{f}$-gradient $\bar{w}$ such that its restriction to $W_{a_1,a_2}$ is an
almost transverse $\bar{f}|$-gradient by Lemma \ref{pajlem}. Furthermore the same is true for the
weak $\omega$-gradients $v_t=tw+(1-t)v$ for $t\in[0,1]$.\\[0.2cm]
Now choose a negative integer $k$ such that $(k+1)b<R$. We set $W=\bar{f}^{-1}[kb,0]$. We claim
that any weak $\omega$-gradient $w$ that lifts to an almost transverse weak $\bar{f}|_W$-gradient satisfies
$b_\omega(-w)<R$. Assume not, then there exists a broken closed orbit $\delta$ of $-w$ with
$\xi(\{\delta\})\geq R$. Let $p\in M$ be a critical point in the image of $\delta$. Lift $p$ to
$\bar{p}\in\bar{f}^{-1}[-b,0]$. The loop $\delta$ lifts to a path $\bar{\delta}$ in $\bar{M}$
starting at $\bar{p}$ and $\bar{\delta}$ is a concatenation of trajectories of $-\bar{w}$ between
critical points and it ends in a translate of $\bar{p}$. But since $\xi(\{\delta\})\geq R>(k+1)b$,
the path $\bar{\delta}$ actually is in $W$. Now $\bar{\delta}$ contradicts almost transversality of
$\bar{w}$ on $W$.\\[0.2cm]
Therefore the family $v_t$ is $R$-controlled and the statement follows by Proposition
\ref{contzeta}.\\[0.2cm]
So now assume $\omega$ is irrational. By Lemma \ref{approx} there is a rational form $\omega'$
agreeing with $\omega$ in a neighborhood of the critical points such that $v$ is an $\omega'$-gradient.
We can also choose $\omega'$ arbitrary close to $\omega$. In particular we can choose $\omega'$
so close that we have \\$\omega_x(v(x))>\frac{1}{2}\omega'_x(v(x))>0$, if $x\in M$ is not in a
neighborhood of a critical point. By compactness and continuity we can also get this for tangent
vectors near $v(x)$. To be more precise, a Riemannian metric on $M$ induces a norm $\|\cdot\|_x$
on every $T_xM$ and we can find an $\varepsilon>0$ such that $\omega_x(X)\geq\frac{1}{2}\omega'_x
(X)>0$ for $X\in T_xM$ with $\|X-v(x)\|_x<\varepsilon$ and all $x$ outside the neighborhood of the
critical points where $\omega$ and $\omega'$ agree.\\[0.2cm]
Now the first part of the proof applies to $\omega'$. We choose a neighborhood $U(v)$ of $v$ such
that every $w\in U(v)$ satisfies
\begin{enumerate}
\item $w$ is an $\omega$-gradient.
\item $\|v(x)-w(x)\|_x<\varepsilon$ for all $x\in M$.
\item $v_t=tw+(1-t)v$ is $(2R)$-controlled with respect to $\omega'$.
\end{enumerate}
It remains to prove that $v_t$ is $R$-controlled with respect to $\omega'$.\\[0.2cm]
So let $\delta$ be a broken closed orbit of some $v_t$. Write $\xi'$ for the homomorphism induced
by $\omega'$. Then
\begin{eqnarray*}
\xi(\{\delta\})&=&\int\limits_\delta\omega\,=\,\sum_{i=1}^k\int\limits_{\gamma_i}\omega\,
=\,-\sum_{i=1}^k\int\limits_{-\infty}^\infty\omega_{\gamma_i(s)}(v_t(\gamma_i(s)))\,ds\\
&\leq&-\sum_{i=1}^k\frac{1}{2}\int\limits_{-\infty}^\infty\omega'_{\gamma_i(s)}
(v_t(\gamma_i(s)))\,ds\\
&=&\frac{1}{2}\int\limits_\delta\omega'\,=\,\frac{1}{2}\xi'(\{\delta\})\,\leq\,\frac{1}{2}\,2R\,=\, R
\end{eqnarray*}
since $b_{\omega'}(-v_t)\leq 2R$. Therefore $(v_t)_{t\in[0,1]}$ is $R$-controlled as a one
parameter family of weak $\omega$-gradients and the statement follows again by Proposition
\ref{contzeta}.
\end{proof}
Corollary \ref{constnc} states that $\zeta$ is constant if there are no critical points, but in
general $\zeta$ is nonconstant, see \cite[Rm.5.4]{schuet}.
\section{Chain homotopy equivalences between Novikov complexes}\label{secwex}
Given two transverse $\omega$-gradients $v,\,w$ for the Morse form $\omega$ we want to describe a
chain homotopy equivalence $\psi_{w,v}$ between the Novikov complexes $C_\ast(\omega,w)$ and
$C_\ast(\omega,v)$ such that the diagram
\begin{equation}\label{chncom}
\begin{array}{rcl}
\multicolumn{3}{c}{\widehat{\mathbb{Z}G}_\xi\otimes_{\mathbb{Z}G}C^\Delta_\ast(\tilde{M})}\\[0.2cm]
\varphi(w)\swarrow& &\searrow\varphi(v)\\[0.2cm]
C_\ast(\omega,w)\hspace{0.5cm}&\stackrel{\psi_{w,v}}{\longrightarrow}&\hspace{0.5cm}
C_\ast(\omega,v)\\[0.2cm] \end{array}
\end{equation}
commutes up to chain homotopy. Then $\tau(\psi_{w,v})=\tau(\varphi(v))-\tau(\varphi(w))$.
Constructions of such equivalences are given in various places in the literature, e.g.\ Latour
\cite[\S 2.21]{latour} gives a description of the torsion which is particularly useful in trying
to show that the torsion of $\varphi(v)$ depends continuously on $v$. In order to show that
(\ref{chncom}) commutes up to chain homotopy using the equivalence of Latour we will give full
proofs for some results in Latour \cite{latour} instead of just refering to \cite{latour} to make
the proof easier to read. Notice that in \cite{schuet} we only showed that $\varphi(v)$ is a chain
homotopy equivalence for $v\in\mathcal{GA}(\omega)$. The fact that $\varphi(v)$ is a chain homotopy
equivalence in general will follow from the fact that $\psi_{v,w}$ is an equivalence once we show
that (\ref{chncom}) commutes up to chain homotopy. A more direct proof can be obtained using the
methods of Pajitnov \cite{pajito}. In fact these methods simplify since the diagram corresponding
to \cite[Diag.(4.1)]{pajito} commutes ``on the nose'' and not just up to chain homotopy.\\[0.2cm]
Let us recall some definitions of Pajitnov \cite{pajiov}. Let $f:W\to[a,b]$ be a Morse function
on a Riemannian cobordism $(W;M_0,M_1)$ and $v$ a transverse $f$-gradient (the Riemannian metric
is only needed to get a metric on the cobordism, but not to specify $v$). If $p$ is a critical point
and $\delta>0$, let $B_\delta(p)$, resp. $D_\delta(p)$ be the image of the Euclidean open, resp.
closed, ball of radius $\delta$ under the exponential map. Here $\delta$ is understood to be so small
that $\exp$ restricts to a diffeomorphism of these balls and so that for different critical points
$p,\,q$ we get $D_\delta(p)\cap D_\delta(q)=\emptyset$.\\[0.2cm]
If $\Phi$ denotes the flow of $v$, we set
\begin{eqnarray*}
B_\delta(p,v)&=&\{x\in W\,|\,\exists t\geq 0\hspace{0.4cm}\Phi(x,t)\in B_\delta(p)\}\\
D_\delta(p,v)&=&\{x\in W\,|\,\exists t\geq 0\hspace{0.4cm}\Phi(x,t)\in D_\delta(p)\}
\end{eqnarray*}
We also define for $i=-1,\ldots\hspace{-1pt},n$
\begin{eqnarray*}
D^i_\delta(v)&=&\bigcup_{{\rm ind}\,p\leq i}D_\delta(p,v)\cup M_0\\
C^i_\delta(v)&=&W-\bigcup_{{\rm ind}\,p\geq i+1} B_\delta(p,-v)\\
C^i(v)&=&W-\bigcup_{{\rm ind}\,p\geq i+1} W^u(p,v)
\end{eqnarray*}
Using a self-indexing Morse function $\phi$ adjusted to $(f,v)$, i.e.\ $v$ is a $\phi$-gradient
and $df=d\phi$ near the critical points and $\partial W$, we get another filtration $W^i=\phi^{-1}
([-\frac{1}{2},i+\frac{1}{2}])$, the one used in Milnor \cite{milnhc}.
\begin{lem}
For $\delta>0$ sufficiently small and $0<\delta_0<\delta$ we have for all $i=-1,\ldots\hspace{-1pt}
,n$
\[D_{\delta_0}^i(v)\subset D^i_\delta(v)\subset W^i\subset C_\delta^i(v)\subset C_{\delta_0}^i(v)\]
and all inclusions are homotopy equivalences.
\end{lem}
\begin{proof}
We can use the flow of $-v$ to define homotopy inverses to the inclusions. For details see
Pajitnov \cite[Prop.2.42]{pajiov}
\end{proof}
Therefore we can use any of these filtrations for the Morse-Smale complex
$C^{MS}_\ast(\tilde{W},\tilde{M}_0;v)$. Since $H_\ast(\tilde{C}^i(v),\tilde{C}^{i-1}(v))$ is the
direct limit of $H_\ast(\tilde{C}^i_\delta(v),\tilde{C}^{i-1}_\delta(v))$ for $\delta>0$ we can
also use $C^i(v)$.\\[0.2cm]
Now if $\Delta$ is a smooth triangulation adjusted to $v$, the chain homotopy equivalence \\
$\varphi(v):C_\ast^\Delta(\tilde{W},\tilde{M}_0)\to C^{MS}_\ast(\tilde{W},\tilde{M}_0;v)$ is just
induced by the inclusion \\$(W^{(i)},W^{(i-1)})\subset (C^i(v),C^{i-1}(v))$. Here $W^{(i)}$ is the
$i$-skeleton of the triangulation.\\[0.2cm]
Let $w$ be another transverse $f$-gradient. To define $\psi_{w,v}:C^{MS}_\ast(\tilde{W},
\tilde{M}_0;w)\to C^{MS}_\ast(\tilde{W},\tilde{M}_0;v)$ let $\Phi:W\to W$ be isotopic to the
identity such that $\Phi(W^u(p,v))\,|\hspace{-7pt}\cap W^s(q,w)$ for ${\rm ind}\,p\geq{\rm ind}\,q$.
The existence is given in Latour \cite[Lm.2.20]{latour}. Furthermore $\Phi$ can be chosen as close
as we like to the identity. Notice that for ${\rm ind}\,q < {\rm ind}\, p$ the intersection is
empty and by compactness we can find a $\delta>0$ such that $D^i_\delta(w)\subset\Phi(C^i(v))$. Let
\[\psi_{w,v}=\tilde{\Phi}^{-1}_\ast\circ j^{-1}_\ast:H_i(\tilde{C}^i(w),\tilde{C}^{i-1}(w))
\stackrel{\simeq}{\longrightarrow}H_i(\tilde{D}^i_\delta(w),\tilde{D}^{i-1}_\delta(w))\longrightarrow
H_i(\tilde{C}^i(v),\tilde{C}^{i-1}(v)).\]
Clearly $\psi_{w,v}$ is a chain map. The $f$-gradients $w,\,v$ have the same critical points, so
we can choose compatible orientations of the stable manifolds.
For ${\rm ind}\,q={\rm ind}\,p$ we have $\Phi(W^u(p,v))\cap W^s(q,w)$ is a finite set and $\psi_{w,v}$
can be expressed by intersection numbers which we denote as $[q:p]\in\mathbb{Z}G$. In particular
we get $[q:q]=1$ and for $f(p)\geq f(q)$ with $p\not=q$ we get $[q:p]=0$, since the intersection
is empty. Thus each $\psi_{w,v}$ can be expressed by an elementary matrix, so $\psi_{w,v}$ is a
simple isomorphism.\\[0.2cm]
Now let $M$ be a closed connected smooth manifold, $\omega$ a Morse form and $v,\,w$ transverse
$\omega$-gradients. Let $\Phi:M\to M$ be isotopic to the identity such that $\Phi(W^u(p,v))
\,\,|\hspace{-8pt}\cap W^s(q,w)$ for ${\rm ind}\,p\geq{\rm ind}\,q$, see Latour \cite[Lm.2.20]{latour}.
Again $\Phi$ can be arbitrarily close to the identity. Choose liftings of the critical points
in the universal cover and orientations of the stable manifolds of $v$. This gives a basis of
$C_\ast(\omega,v)$ and we can choose a corresponding basis for $C_\ast(\omega,w)$.
\begin{prop}\label{latequ}
If $p,\,q$ are critical points of the same index, the intersection number \\$[q:p]\in\widehat{\mathbb{Z}G}_\xi$
is well defined and $\psi_{w,v}:C_\ast(\omega,w)\to C_\ast(\omega,v)$ given by
$\psi_{w,v}(q)=\sum[q:p]\,p$ is an isomorphism of chain complexes with $\tau(\psi_{w,v})\in
\overline{W}$.
\end{prop}
\begin{proof}
Assume $\omega$ is rational. Let $p^\ast\omega=df$ with $f:\bar{M}\to\mathbb{R}$ having $0\in\mathbb{R}$
as a regular value. Also let $\tilde{f}:\tilde{M}\to\mathbb{R}$ be the composition with the universal
covering projection. We can assume that the $b$ from Section \ref{novcom} is 1 and the liftings
of the critical points are chosen in $\tilde{f}^{-1}([-1,0])$.\\[0.2cm]
Define
\[C_m^i(v)=\bar{\Phi}\left(f^{-1}((-\infty,0])-\bigcup_{{\rm ind}\,r\geq i+1 \atop f(r)\geq -m}
W^u(r,\bar{v})\right)\]
where $\bar{v}$ is the lift of $v$ to $\bar{M}$ and the same for $\bar{\Phi}$. Also let
\[D^i_{m,\delta_m}(w)=f^{-1}((-\infty,-m]\cup\bigcup_{{\rm ind}\,r\leq i}D_\delta(r,\bar{w})\]
with $\delta_m>0$ so small that $D^i_{m,\delta_m}(w)\subset C_m^i(v)$.
Now define \\$C^{MS}_i(m,w)=H_i(\tilde{D}^i_{m,\delta_m}(w),\tilde{D}^{i-1}_{m,\delta_m}(w))$ and
$C^{MS}_i(m,v)=H_i(\tilde{C}^i_m(v),\tilde{C}^{i-1}_m(v))$. Both complexes calculate the homology
of $(\tilde{f}^{-1}([-m,0]),\tilde{f}^{-1}(\{-m\}))$. The exact case described above gives a
chain isomorphism $\psi^m_{w,v}:C_i^{MS}(m,w)\to C_i^{MS}(m,v)$ such that the diagram
\[\begin{array}{ccc}
C_\ast^{MS}(m,w)&\longleftarrow&C_\ast^{MS}(m+1,w)\\[0.3cm]
\Big\downarrow\psi^m_{w,v}& &\Big\downarrow\psi^{m+1}_{w,v}\\[0.3cm]
C_\ast^{MS}(m,v)&\longleftarrow&C_\ast^{MS}(m+1,v)\end{array}\]
commutes. In fact all the arrows are just induced by inclusion. Passing to the inverse limit gives
almost the Novikov complex; we only look at $\tilde{f}^{-1}((-\infty,0])\subset\tilde{M}$.
But $\lim\limits_{\longleftarrow} C_\ast^{MS}(m,w)$ is a finitely generated free
$\widehat{\mathbb{Z}G}_\xi^0$ complex\footnote{Notice that we considered $C_\ast^{MS}(m,w)$ as a
$\mathbb{Z}$ complex, so the inverse limit is also just a $\mathbb{Z}$ complex, but it carries
extra structure as a free $\widehat{\mathbb{Z}G}^0_\xi$ complex generated by the critical points
of $\omega$.}, where $\widehat{\mathbb{Z}G}^0_\xi$ is the subring of $\widehat{\mathbb{Z}G}_\xi$
consisting of elements $\lambda$ with $\|\lambda\|\leq 1$.\\[0.2cm]
Now $C_\ast(\omega,w)=\widehat{\mathbb{Z}G}_\xi\otimes_{\widehat{\mathbb{Z}G}^0_\xi}
\lim\limits_{\longleftarrow} C_\ast^{MS}(m,w)$ and similarly for $C_\ast(\omega,v)$. The chain map
$\psi_{w,v}={\rm id}_{\widehat{\mathbb{Z}G}_\xi}\otimes_{\widehat{\mathbb{Z}G}^0_\xi}\lim\limits
_{\longleftarrow}\psi^m_{w,v}$ is represented by intersection numbers, since the $\psi^m_{w,v}$ are.
In particular we have $[q:p]\in\widehat{\mathbb{Z}G}_\xi$. Also $[q:q]=1-a_q$ with $\|a_q\|<1$
and $[q:p]=b_{qp}-a_{qp}$ with $\|a_{qp}\|<1$ and $b_{qp}$ is the coefficient of $\psi^1_{w,v}$.
So we can order the critical points such that the matrix of $\psi_{w,v}$ is of the form $I-O-A$,
where $O$ is nilpotent and $A$ satisfies $\|A_{ij}\|<1$ for all entries. The matrix $I+O+O^2+\ldots$
is elementary and $(I-O-A)\cdot(I+O+O^2+\ldots)=I-A'$ where the entries of $A'$ satisfy $\|A'_{ij}
\|<1$. Therefore $\psi_{v,w}$ is an isomorphism of chain complexes and $\tau(\psi_{w,v})\in
\overline{W}$.\\[0.2cm]
It remains to prove the proposition for irrational $\omega$. We can assume that there exists a
rational approximation $\omega'$ that agrees with $\omega$ near the critical points such that $v$
and $w$ are also $\omega'$-gradients for we can find a sequence $w=w_0,w_1,\ldots\hspace{-1pt},w_k
=v$ of $\omega$-gradients such that $w_i$ and $w_{i+1}$ have a common rational approximation.\\[0.2cm]
Let $\xi':G\to\mathbb{R}$ be the homomorphism induced by $\omega'$. By the rational case above
we find a chain isomorphism $\psi_{w,v}'$ for the $\widehat{\mathbb{Z}G}_{\xi'}$ Novikov complexes
$C_\ast(\omega',w)$ and $C_\ast(\omega',v)$. If we can show that the matrix entries of $\psi_{w,v}'$
lie in $\widehat{\mathbb{Z}G}_{\xi'}\cap\widehat{\mathbb{Z}G}_\xi$, then $\psi_{w,v}'$ induces a
chain map $\psi_{w,v}:C_\ast(\omega,w)\to C_\ast(\omega,v)$ by the remarks at the end of Section
\ref{novcom}. Notice that the entries are intersection numbers $[q:p]$. So $[q:p](g)\not=0$ gives
a point $\tilde{x}\in\tilde{M}$, a trajectory $\gamma_1$ of $-\tilde{w}$ from $\tilde{q}$ to
$\tilde{x}$ and a trajectory $\gamma_2$ of $-\tilde{v}$ from $\tilde{\Phi}^{-1}(\tilde{x})$ to
$\tilde{\Phi}^{-1}(g\tilde{p})$. Now $[q:p]\in\widehat{\mathbb{Z}G}_{\xi'}\cap\widehat{\mathbb{Z}G}
_\xi$ follows from the next lemma.\\[0.2cm]
To see that $\psi_{w,v}$ is an isomorphism with $\tau(\psi_{w,v})\in\overline{W}$ notice that in
the irrational case we can choose liftings of the critical points of $\omega$ in $\tilde{f}^{-1}(I_
\varepsilon)$, where $I_\varepsilon\subset\mathbb{R}$ is an arbitrarily small interval. Then the
matrix of $\psi_{w,v}$ in a basis corresponding to these critical points is of the form $I-A$
with $\|A_{ij}\|<1$ for all the entries of $A$, compare Latour \cite[\S 2.23]{latour}
\end{proof}
\begin{lem}
Let $\omega_1,\,\omega_2$ be Morse forms  that agree near the common set of critical points
with corresponding homomorphisms $\xi_1,\xi_2:G\to\mathbb{R}$. Let $v,\,w$ be both $\omega_1$- and
$\omega_2$-gradients. Then there exist constants $A,\,B\in\mathbb{R}$ with $A>0$ such that whenever
for $g\in G$ there exist critical points $p,\,q$, a point $\tilde{x}\in\tilde{M}$, a trajectory
$\tilde{\gamma}_1$ of $-\tilde{w}$ from $\tilde{q}$ to $\tilde{x}$ and a trajectory
$\tilde{\gamma}_2$ of $-\tilde{v}$
from $\tilde{\Phi}^{-1}(\tilde{x})$ to $\tilde{\Phi}^{-1}(g\tilde{p})$, then $\xi_1(g)\leq A\xi_2
(g)+B$.
\end{lem}
\begin{proof}
For every pair of critical points $p,\,q$ of $\omega_1$ we can choose a path $\tilde{\gamma}_{pq}$
in$\tilde{M}$ from $\tilde{p}$ to $\tilde{q}$. Then there is a constant $K>0$ such that
$|\int_{\gamma_{pq}}\omega_i|\leq K$ for $i=1,2$ and all pairs of critical points. Let $\Theta:
M\times I\to M$ be the isotopy between id and $\Phi$. For every $y\in M$ we get a path $\gamma_y(t)
=\Theta(y,t)$ from $y$ to $\Phi(y)$. By compactness we can also assume $|\int_{\gamma_y}\omega_i|
\leq K$ for $i=1,2$ and all $y\in M$. Since $\omega_1$ and $\omega_2$ agree near the critical points
there exists a $C\in(0,1)$ such that
\[\omega_1(v(x))\geq C\omega_2(v(x))\mbox{ and }\omega_1(w(x))\geq C\omega_2(w(x))\mbox{ for all }
x\in M\]
again by compactness. Now let $g\in G$ be as in the statement. Then
\begin{eqnarray*}
\xi_2(g)&=&\int_{\gamma_{qp}}\omega_2+\int_{\gamma_1}\omega_2+\int_{\gamma_x}\omega_2+\int_{\gamma_2}\omega_2\\
&\geq &-2K-\int_{-\infty}^{b_1}(\omega_2)_{\gamma_1(t)}(w(\gamma_1(t)))\,dt-
\int_{a_1}^\infty(\omega_2)_{\gamma_2(t)}(v(\gamma_2(t)))\,dt\\
&\geq&-2K-C\left(\int_{-\infty}^{b_1}(\omega_1)_{\gamma_1(t)}(w(\gamma_1(t)))\,dt+
\int_{a_1}^\infty(\omega_1)_{\gamma_2(t)}(v(\gamma_2(t)))\,dt\right)\\
&\geq&-2K-2KC+C\left(\int_{\gamma_{qp}}\omega_1+\int_{\gamma_1}\omega_1+\int_{\gamma_x}\omega_1+
\int_{\gamma_2}\omega_1\right)\\
&=&-2K(1+C)+C\xi_1(g)
\end{eqnarray*}
which gives the result.
\end{proof}
To show that (\ref{chncom}) commutes up to chain homotopy let us start with the exact case again,
i.e.\ we have a compact cobordism and a Morse function $f:W\to[a,b]$. We use the same notation as
before. Let $\Delta$ be a smooth triangulation adjusted to $w$ and $\Phi_\ast v=d\Phi^{-1}\circ v
\circ \Phi$, this is possible by \cite[\S A.1]{schuet}. So for every $k$-simplex $\sigma$ we have
$\sigma\,\,|\hspace{-8.5pt}\cap W^u(p,w)$ and \\$\sigma\,\,|\hspace{-8.5pt}\cap\Phi(W^u(p,v))$ if ${\rm ind}
\,p\geq k$.
\begin{prop}\label{exacom}
The chain maps $\psi_{w,v}\circ\varphi(w)$ and $\varphi(v)$ are chain homotopic.
\end{prop}
\begin{proof}
Let $\Theta_w:\tilde{W}\times\mathbb{R}\to\tilde{W}$ be induced by the flow of $-w$, i.e.\ stop
once the boundary is reached. There is a $\delta>0$ such that $D_\delta^k(w)\subset\Phi(C^k(v))$.
Since $\Delta$ is adjusted to $w$ there is a $K>0$ such that $\Theta_w(\tilde{W}^{(k)},K)\subset
\tilde{D}^k_\delta(w)$, where $W^{(k)}$ is the $k$-skeleton of the triangulation. Furthermore
$\Theta_w$ gives a homotopy between id and $\Theta(\cdot,K)$.\\[0.2cm]
Since $\Delta$ is adjusted to $\Phi_\ast v$ we have $W^{(k)}\subset\Phi(C^k(v))$. We can modify the
homotopy away from the endpoints to get a homotopy $h:\tilde{W}\times I\to\tilde{W}$ between
id and $\Theta_w(\cdot,K)$ such that $h(\tilde{W}^{(k)}\times I)\subset\Phi(C^{k+1}(v))$. The
modifications are done skeleton by skeleton, compare the proof of \cite[Lm.A.2]{schuet} and can
be done arbitrarily close to the original homotopy. Now define $H:C_k^\Delta(\tilde{W},\tilde{M}_0)
\to C^{MS}_{k+1}(\tilde{W},\tilde{M}_0;v)$ be sending $\tilde{\sigma}$ to $(-1)^k\tilde{\Phi}_\ast
^{-1}h_\ast[\tilde{\sigma}\times I]\in H_{k+1}(\tilde{C}^{k+1}(v),\tilde{C}^k(v))$. Then
\begin{eqnarray*}
\partial H+H\partial(\tilde{\sigma})&=&(-1)^k \tilde{\Phi}^{-1}_\ast h_\ast\partial[\tilde{\sigma}
\times I]+(-1)^{k-1}\tilde{\Phi}^{-1}_\ast h_\ast[\partial\tilde{\sigma}\times I]\\
&=&\tilde{\Phi}_\ast^{-1}h_\ast[\tilde{\sigma}\times 1]-\tilde{\Phi}_\ast^{-1}h_\ast[\tilde{\sigma}
\times 0]\\
&=&\tilde{\Phi}^{-1}_\ast\Theta_{w\ast}[\tilde{\sigma}\times K]-\tilde{\Phi}_\ast^{-1}[\tilde{\sigma}]\\
&=&\psi_{w,v}\circ\varphi(w)(\tilde{\sigma})-\varphi(v)(\tilde{\sigma}).
\end{eqnarray*}
Notice that $\Theta_{w\ast}[\tilde{\sigma}\times K]\in H_k(\tilde{D}^k_\delta(w),\tilde{D}^{k-1}
_\delta(w))$ represents $\varphi(w)(\tilde{\sigma})$ and using $\tilde{\Phi}_\ast^{-1}$ gives
$\psi_{w,v}$.
\end{proof}
\begin{prop}\label{latcom}
Diagram (\ref{chncom}) commutes up to chain homotopy.
\end{prop}
\begin{proof}
Assume $\omega$ is rational. We use the notation from the proof of Proposition \ref{latequ}. We
can assume that $\Delta$ contains $f^{-1}(\{0\})$ as a subcomplex. Let us also set $\tilde{M}_m
=\tilde{f}^{-1}([-m,0])$.\\[0.2cm]
Let $H^m:C^\Delta_k(\tilde{M}_m,\tilde{f}^{-1}(\{-m\}))\to C^{MS}_{k+1}(m,v)$ be the chain homotopy
given by Proposition \ref{exacom}. Actually in the nonexact case it comes from a homotopy $h_m:
\tilde{f}^{-1}((-\infty,0])\times I\to\tilde{f}^{-1}((-\infty,0])$ and it satisfies
$h_m(\tilde{\sigma}_k\times I)\subset\tilde{\Phi}(\tilde{C}^{k+1}_m(v))$ and $h_m(\tilde{\sigma}_k
\times 1)\subset\tilde{D}^k_{m,\delta_m}(w)$.\\[0.2cm]
We want to get a chain homotopy $H^{m+1}:C^\Delta_\ast(\tilde{M}_{m+1},\tilde{f}^{-1}(\{-m-1\}))
\to C^{MS}_{k+1}(m+1,v)$ based on $H^m$. First we need $h_{m+1}(\tilde{\sigma}_k\times 1)\subset
\tilde{D}^k_{m+1,\delta_{m+1}}(w)$. Notice that $\delta_{m+1}\leq \delta_m$. So we take the
homotopy $h_m$ and flow along $-\tilde{w}$ for a little bit longer. Call this homotopy $h_{m+1}'$.
Then $h_{m+1}'(\tilde{\sigma}_k\times I)\subset \tilde{\Phi}(\tilde{C}^{k+1}_m(v))$, but not
necessarily $\subset\tilde{\Phi}(\tilde{C}^{k+1}_{m+1}(v))$. We need to adjust the homotopy slightly
to achieve this. So do this skeleton by skeleton to get a homotopy $h_{m+1}$ so close to $h_{m+1}'$
that passing $\tilde{\sigma}_k\times I$ from $h_{m+1}'$ to $h_{m+1}$ is done within $\tilde{\Phi}
(\tilde{C}^{k+1}_m(v))$.\\[0.2cm]
Then if $H^{m+1}$ is induced by $h_{m+1}$ as in the proof of Proposition \ref{exacom} we get the
commutative diagram
\[\begin{array}{ccc}
C^\Delta_k(\tilde{M}_m,\tilde{f}^{-1}(\{-m\}))&\longleftarrow& C_k^{\Delta}(\tilde{M}_{m+1},
\tilde{f}^{-1}(\{-m-1\}))\\[0.3cm]
\Big\downarrow H^m& &\Big\downarrow H^{m+1}\\[0.3cm]
C^{MS}_{k+1}(m,v)&\longleftarrow&C_{k+1}^{MS}(m+1,v)
\end{array}\]
Passing to the inverse limit as in the proof of Proposition \ref{latequ} gives the result in
the rational case.\\[0.2cm]
For the irrational case notice that nonzero terms of the chain homotopy give a trajectory of
$-w$ from some $x\in\sigma_k$ to a $y\in M$ and a trajectory of $-v$ from $\Phi^{-1}(y)$ to a
critical point. Thus we can use a similar approximation argument as in the proof of Proposition
\ref{latequ}, we omit the details.
\end{proof}
\section{Proof of the main theorem, Part 2}
In \cite{schuet} we have shown that $\varphi(v)$ is a chain homotopy equivalence for $v\in
\mathcal{GA}(\omega)$ and that $\tau(\varphi(v))\in\overline{W}$. By Proposition \ref{latequ}
and Proposition \ref{latcom} we now get that this also holds for any transverse $\omega$-gradient
$v$. But we want to show that the torsion actually depends continuously on the vector field. To
do this let us first put a topology on $\overline{W}$. Denote by $U$ the subgroup of units of
$\widehat{\mathbb{Z}G}_\xi$ that consists of elements $1-a$ with $\|a\|<1$. As a subset of
$\widehat{\mathbb{Z}G}_\xi$ it carries a natural topology. Now $U$ surjects onto $\overline{W}$
so we give $\overline{W}$ the quotient topology. Notice that both $U$ and $\overline{W}$ are
topological groups.\\[0.2cm]
For a Morse form $\omega$ we let $\mathcal{G}_t(\omega)$ be the space of transverse $\omega$-gradients
with the $C^0$-topology.
\begin{theo}\label{torcont}
Let $\omega$ be a Morse form on the closed connected smooth manifold $M$. Then the map $\mathcal{T}
:\mathcal{G}_t(\omega)\to\overline{W}$ given by $\mathcal{T}(v)=\tau(\varphi(v))$ is continuous.
\end{theo}
\begin{proof} For $R<0$ let $U_R=\{1-a\,|\,\|a\|< \exp R\}$. The collection $(U_R)_{R<0}$ forms
a neighborhood basis of $1\in U$, so $(\tau(U_R))_{R<0}$ forms a neighborhood basis of $0\in
\overline{W}$.\\[0.2cm]
Let $v\in\mathcal{G}_t(\omega)$. To see that $\mathcal{T}$ is continuous, we have to find for every
$R<0$ a neighborhood $\mathcal{U}$ of $v$ such that $w\in\mathcal{U}$ satisfies $\tau(\varphi(v))-
\tau(\varphi(w))\in\tau(U_R)$. By Proposition \ref{latcom} we have $\tau(\varphi(v))-
\tau(\varphi(w))=\tau(\psi_{w,v})$.\\[0.2cm]
Assume that $\omega$ is rational. Let $\bar{M}$ be the infinite cyclic covering space corresponding
to $\ker \xi$ and $f:\bar{M}\to\mathbb{R}$ such that $df$ is the pullback of $\omega$ and $0\in\mathbb{R}$
a regular value. For simplicity assume that the $b$ from Section \ref{novcom} is 1. Since $v$ is
transverse, so is its lift $\bar{v}$ to $\bar{M}$. Choose an integer $m$ with $m+1<R$. We can find
a self-indexing Morse function $\phi:f^{-1}([m,0])\to[-\frac{1}{2},n+\frac{1}{2}]$ such that
$\bar{v}|$ is a $\phi$-gradient and $d\phi=df$ near the critical points and $f^{-1}(\{m,0\})$.
Pajitnov \cite[Lm.2.74]{pajiov} gives a neighborhood $\mathcal{U}$ of $v$ in $\mathcal{G}_t(\omega)$
such that every $w\in\mathcal{U}$ lifts to a $\phi$-gradient on $f^{-1}([m,0])$. Now for
${\rm ind}\,q\leq{\rm ind}\,p$ with $q\not=p$ we get $W^s(q,\bar{w})\cap W^u(p,\bar{v})\cap
f^{-1}([m,0])=\emptyset$, since $\phi(W^s(q,\bar{w})\cap f^{-1}([m,0])-\{q\})\subset[-\frac{1}{2},
{\rm ind}\,q)$ and $\phi(W^u(p,\bar{v})\cap f^{-1}([m,0])-\{p\})\subset({\rm ind}\,p,n+\frac{1}{2}]$.
By choosing the isotopy $\Phi$ of $M$ close enough to the identity we still have $\bar{\Phi}
(W^u(p,\bar{v}))\cap W^s(q,\bar{w})\cap f^{-1}([m,0])=\emptyset$. Now we choose liftings of the
critical points within $f^{-1}([-1,0])$ to get a basis for the Novikov complex. For every $w\in
\mathcal{U}$ the coefficients $[q:p]$ of $\psi_{w,v}$ then have the property that any $g\in G$
with $[q:p](g)\not=0$ implies $\xi(g)<R$, compare the proof of Theorem \ref{zetacont}. Thus
$\tau(\psi_{w,v})$ is represented by a matrix $I-A$ where $\|A_{ij}\|<\exp R$ for all entries of
$A$. By Gau\ss \ elimination we see that $\tau(\psi_{w,v})=\tau(1-a)$ with $\|a\|<\exp R$,
so $\tau(\psi_{w,v})\in\tau(U_R)$ for all $w\in\mathcal{U}$.\\[0.2cm]
The irrational case is now derived from the rational case by analogy to the proof of Theorem
\ref{zetacont}, we omit the details.
\end{proof}
Now we can finally drop the cellularity condition on the vector fields in Theorem \ref{mtheo}
to get
\begin{theo}\label{rmtheo}
Let $\omega$ be a Morse form on a smooth connected closed manifold $M^n$. Let $\xi:G\to\mathbb{R}$
be induced by $\omega$ and let $C^\Delta_\ast(\tilde{M})$ be the simplicial $\mathbb{Z}G$ complex
coming from a smooth triangulation of $M$. For every transverse $\omega$-gradient $v$ there is a
natural chain homotopy equivalence $\varphi(v):\widehat{\mathbb{Z}G}_\xi\otimes_{\mathbb{Z}G}
C_\ast^\Delta(\tilde{M})\to C_\ast(\omega,v)$ given by (\ref{formula}) whose torsion $\tau(\varphi(v))$
lies in $\overline{W}$ and satisfies
\[\mathfrak{DT}(\tau(\varphi(v)))=\zeta(-v).\]
\end{theo}
\begin{proof}
Clearly the homomorphism $\mathfrak{DT}:\overline{W}\to\widehat{H\!H}_1(\mathbb{Z}G)_\xi$ is
continuous. So the statement follows from the theorems \ref{mtheo}, \ref{zetacont} and \ref{torcont},
since $\mathcal{GA}(\omega)$ is dense in $\mathcal{G}_t(\omega)$.
\end{proof}
Let us obtain a commutative version of Theorem \ref{rmtheo}. Instead of the universal covering space
we look at the universal abelian covering space $\overline{M}$. We set $H=H_1(M)$. Then
$C_\ast^\Delta(\overline{M})=\mathbb{Z}H\otimes_{\mathbb{Z}G}C_\ast^\Delta(\tilde{M})$. If we set
$\overline{C}_\ast(\omega,v)=\widehat{\mathbb{Z}H}_{\bar{\xi}}\otimes_{\widehat{\mathbb{Z}G}_\xi}
C_\ast(\omega,v)$, we get the Novikov complex corresponding to $\overline{M}$. Then $\bar{\varphi}(v)
=\mbox{id}\otimes_{\widehat{\mathbb{Z}G}_\xi}\varphi(v):\widehat{\mathbb{Z}H}_{\bar{\xi}}\otimes_{\mathbb{Z}H}
C_\ast^\Delta(\overline{M})\to\overline{C}_\ast(\omega,v)$ is a chain homotopy equivalence. Denote
the subgroup of $K_1(\widehat{\mathbb{Z}H}_{\bar{\xi}})$ consisting of units of the form $1-a$, where
$\|a\|<1$ by $W'$.\\[0.2cm]
To define a commutative zeta function let $\widehat{\mathbb{Q}H}^-_{\bar{\xi}}=\{\lambda\in
\widehat{\mathbb{Q}H}_{\bar{\xi}}\,|\,\|\lambda\|<1\}$, a subgroup of $\widehat{\mathbb{Q}H}_{\bar{\xi}}$.
Notice that $\varepsilon(\eta(-v))\in\widehat{\mathbb{Q}H}^-_{\bar{\xi}}$, where $\varepsilon$ is
the augmentation. We define $\exp:\widehat{\mathbb{Q}H}^-_{\bar{\xi}}\to
1+\widehat{\mathbb{Q}H}^-_{\bar{\xi}}$ by $\exp(\lambda)=\sum\limits_{m=0}^\infty \dfrac{\lambda^m}{m!}$.
\begin{defi}\em
Let $\omega$ be a Morse form and $v$ an $\omega$-gradient with
$b_\omega(-v)=-\infty$. Then we define the \em zeta function \em of $-v$ to be
\[\bar{\zeta}(-v)=\exp\circ\,\varepsilon(\eta(-v))\in
1+\widehat{\mathbb{Q}H}^-_{\bar{\xi}}.\] \end{defi}
Notice that this coincides with the formula for a zeta function given in Fried \cite{fried}.
\begin{coro}\label{cmtheo}
Let $\omega$ be a Morse form and $v$ a transverse $\omega$-gradient. Then there is a natural chain
homotopy equivalence $\bar{\varphi}(v):\widehat{\mathbb{Z}H}_{\bar{\xi}}\otimes_{\mathbb{Z}H}
C_\ast^\Delta(\overline{M})\to\overline{C}(\omega,v)$ whose torsion lies in $W'$ and that satisfies
\[\det(\tau(\bar{\varphi}(v)))=\bar{\zeta}(-v).\]
\end{coro}
\begin{proof}
The composition $\overline{W}\stackrel{\mathfrak{DT}}{\longrightarrow}\widehat{H\!H}_1(\mathbb{Z}G)_\xi
\stackrel{l}{\longrightarrow}\widehat{\mathbb{R}\Gamma}_\xi$ induces the homomorphism \\
$\mathfrak{L}:\overline{W}\to\widehat{\mathbb{Q}\Gamma}_\xi^-$ from \cite[\S 3.2]{schuet}, compare
Section \ref{hhofnov}. By \cite[Prop.3.4]{schuet}, $\mathfrak{L}$ is a power series of a logarithm,
so we get $\exp\circ\,\varepsilon\circ\mathfrak{L}
([1-a])=\det\circ\,\varepsilon_\ast([1-a])$, where $\varepsilon_\ast:\overline{W}\to W'$ is induced
by the augmentation $\widehat{\mathbb{Z}G}_\xi\to\widehat{\mathbb{Z}H}_{\bar{\xi}}$.
By (\ref{etazeta}) and Theorem \ref{rmtheo} we get
\begin{eqnarray*}
\bar{\zeta}(-v)&=&\exp\circ\,\varepsilon\circ l(\zeta(-v))\,\,
=\,\,\exp\circ\,\varepsilon\circ l\circ\mathfrak{DT}(\tau(\varphi(v)))\\
&=&\exp\circ\,\varepsilon\circ\mathfrak{L}(\tau(\varphi(v)))\,\,
=\,\,\det(\tau(\bar{\varphi}(v)))
\end{eqnarray*}
\end{proof}
\section{The zeta function vs.\ the eta function}
In the commutative case the zeta and the eta function carry the same information since we have
\[\bar{\zeta}(-v)=\exp\bar{\eta}(-v)\hspace{1cm}\mbox{and}\hspace{1cm}\bar{\eta}(-v)=\log\bar{\zeta}
(-v).\]
We have seen in Section \ref{novcom} that the noncommutative zeta function determines the
noncommutative eta function via $\eta(-v)=l(\zeta(-v))$. It is natural to ask whether the zeta
function is determined by the eta function as in the commutative case or if it actually carries
more information than the eta function.\\[0.2cm]
Let us define a rational version of the noncommutative zeta function. The ring homomorphism $i:\mathbb{Z}G
\to \mathbb{Q}G$ induces a map on Hochschild homology $i_\ast:H\!H_\ast(\mathbb{Z}G)\to H\!H_\ast
(\mathbb{Q}G)$. Since $\mathbb{Q}$ is a flat $\mathbb{Z}$ module and $\mathbb{Q}\otimes\mathbb{Q}
\simeq \mathbb{Q}$ we see that $H\!H_\ast(\mathbb{Q}G)\simeq \mathbb{Q}\otimes H\!H_\ast(\mathbb{Z}
G)$. Furthermore we get a direct sum decomposition of $C_\ast(\mathbb{Q}G,\mathbb{Q}G)$ as in
Section \ref{hhfg} and we can complete $H\!H_\ast(\mathbb{Q}G)$ to $\widehat{H\!H}_\ast(\mathbb{Q}
G)_\xi$. The homomorphism $i_\ast$ extends to $\hat{\imath}_\ast:\widehat{H\!H}_\ast(\mathbb{Z}G)_\xi
\to \widehat{H\!H}_\ast(\mathbb{Q}G)_\xi$ and we define the \em rational noncommutative zeta
function \em by
\[\zeta_{\mathbb{Q}}(-v)=\hat{\imath}_\ast\zeta(-v)\in\widehat{H\!H}_1(\mathbb{Q}G)_\xi.\]
It is easy to see that $l:\widehat{H\!H}_1(\mathbb{Z}G)_\xi\to\widehat{\mathbb{R}\Gamma}_\xi=
\widehat{H\!H}_0(\mathbb{R}G)_\xi$ factors through $\widehat{H\!H}_1(\mathbb{Q}G)_\xi$ as
$l=l_{\mathbb{Q}}\circ\hat{\imath}_\ast$. For $\gamma\in\Gamma$ define $e_\gamma:C_0(\mathbb{Q}G,
\mathbb{Q}G)_\gamma\to C_1(\mathbb{Q}G,\mathbb{Q}G)_\gamma$ by $e_\gamma(g)=1\otimes g$. By
Lemma \ref{ezlem} this induces a homomorphism $e:\widehat{H\!H}_0(\mathbb{Q}G)_\xi\to
\widehat{H\!H}_1(\mathbb{Q}G)_\xi$ with $l_\mathbb{Q}\circ e(x)=x$ for $x\in \widehat{\mathbb{Q}
\Gamma}_\xi^-$. Notice that $\eta(-v)\in\widehat{\mathbb{Q}\Gamma}_\xi^-$.
\begin{prop}
Let $\omega$ be a Morse form on the closed connected smooth manifold $M$ and $v$ an almost
transverse $\omega$-gradient. Then $\zeta_\mathbb{Q}(-v)=e(\eta(-v))$.
\end{prop}
\begin{proof}
Assume that the closed orbits of $v$ are nondegenerate. Then for a closed orbit $\gamma$ of $-v$
of multiplicity $m$, we get a summand $\frac{\varepsilon(\gamma)}{m}\{\gamma\}$ in $\eta(-v)$.
Let $g\in G$ be so that $g^m$ represents the conjugacy class $\{\gamma\}$. Then $e(\frac{\varepsilon
(\gamma)}{m}\{\gamma\})=[\frac{\varepsilon(\gamma)}{m}\,1\otimes g^m]$, but $\frac{1}{m}\otimes
g^m$ is homologous to $g^{m-1}\otimes g$ by Lemma \ref{homtrace}. Therefore $e(\frac{\varepsilon
(\gamma)}{m}\{\gamma\})=\varepsilon(\gamma)I(\gamma)$ (recall the Nielsen-Fuller series from
Section \ref{gradflws}) and we get the result.\\[0.2cm]
The general case now follows by continuity, compare the end of Section \ref{novcom}.
\end{proof}
To simplify notation let $H_\gamma=H_1(C_\ast(\mathbb{Z}G,\mathbb{Z}G)_\gamma)$ for
$\gamma\in\Gamma$. Projection gives homomorphisms $p_\gamma:\widehat{H\!H}_1(\mathbb{Z}G)
_\xi\to H_\gamma$ and $p_{\mathbb{Q},\gamma}:\widehat{H\!H}_1(\mathbb{Q}G)
_\xi\to \mathbb{Q}\otimes H_\gamma$. It follows from Section \ref{gradflws} that $p_\gamma(\zeta(
-v))$ is generated by homology classes of the form $[g^{k-1}\otimes g]$ where $k\geq 1$ and
$\gamma(g^k)=\gamma$. It is possible that $p_\gamma(\zeta(-v))$ is a torsion element, so that
$p_{\mathbb{Q},\gamma}(\zeta_\mathbb{Q}(-v))=0$. But to produce torsion we need closed orbits of
multiplicity $>1$. It is not clear to the author whether at the $\mathbb{Z}$-level the zeta function
carries more information than the eta function.


\begin{thebibliography}{99}
\bibitem{abrah}R. Abraham and J. Robbin, Transversal mappings and flows, W.A. Benjamin, 1967.
\bibitem{brown}R. Brown, The Lefschetz fixed point theorem, Scott Foresman, 1971.
\bibitem{eisenb}D. Eisenbud, Commutative Algebra with a view toward algebraic geometry, Springer, 1994.
\bibitem{farran}M. Farber and A. Ranicki, The Morse-Novikov theory of circle-valued functions and noncommutative localization, Tr. Mat. Inst. Steklova 225 (1999) 381-388.
\bibitem{fried}D. Fried, Homological identities for closed orbits, Inv. Math. 71 (1983), 419-442.
\bibitem{fuller}F. Fuller, An index of fixed point type for periodic orbits, Amer. J. Math. 89 (1967), 133-148.
\bibitem{geonia}R. Geoghegan and A. Nicas, Parametrized Lefschetz-Nielsen fixed point theory and Hochschild homology traces, Amer. J. Math. 116 (1994), 397-446.
\bibitem{geonic}R. Geoghegan and A. Nicas, Trace and torsion in the theory of flows, Topology 33 (1994), 683-719.
\bibitem{genisc}R. Geoghegan, A. Nicas and D. Sch\"utz, Obstructions to homotopy invariance in parametrized fixed point theory, in: Geometry and Topology: Aarhus, Editors K. Grove, I. Madsen and E. Pedersen, Contemp. Math. 258 (2000), 157-175.
\bibitem{hutcth}M. Hutchings, Reidemeister torsion in generalized Morse theory, Harvard University Ph.D. thesis, 1998.
\bibitem{hutchi}M. Hutchings, Reidemeister torsion in generalized Morse theory, available as math.DG/9907066, to appear in Forum Math.
\bibitem{hutlee}M. Hutchings and Y-J. Lee, Circle-valued Morse theory, Reidemeister torsion, and Seiberg-Witten invariants of three manifolds, Topology 38 (1999), 861-888.
\bibitem{hutle2}M. Hutchings and Y-J. Lee, Circle-valued Morse theory and Reidemeister torsion, Geom. Topol. 3 (1999), 369-396.
\bibitem{igusa}K. Igusa, What happens to Hatcher and Wagoner's formula for $\pi_0C(M)$ when the first Postnikov invariant is nontrivial?, Algebraic K-theory, Number theory, Geometry and Analysis, Lecture notes in Math. vol 1046, Springer 1984, 104-172.
\bibitem{latour}F. Latour, Existence de 1-formes ferm\'ees non singuli\`eres dans une classe de cohomologie de de Rham, Publ. IHES No.80 (1994), 135-194.
\bibitem{milnhc}J. Milnor, Lectures on the h-cobordism theorem, Princeton University Press, 1965.
\bibitem{pajito}A. Pazhitnov, On the Novikov complex for rational Morse forms, Ann. Fac. Sci. Toulouse 4 (1995), 297-338.
\bibitem{pajisp}A. Pajitnov, Incidence coefficients in the Novikov complex for Morse forms: rationality and exponential growth properties, available as math.DG/9604004.
\bibitem{pajirn}A. Pajitnov, Simple homotopy type of the Novikov complex and Lefschetz $\zeta$-functions of the gradient flow, Russ. Math. Surveys 54 (1999), 119-169.
\bibitem{pajiov}A. Pajitnov, $C^0$-generic properties of boundary operators in the Novikov complex, Pseudoperiodic topology, Amer. Math. Soc. Transl. Ser. 2, 197 (1999), 29-115.
\bibitem{pajitn}A. Pajitnov, Closed orbits of gradient flows and logarithms of non-abelian Witt vectors, K-Theory 21 (2000), 301-324.
\bibitem{pajran}A. Pajitnov and A. Ranicki, The Whitehead group of the Novikov ring, K-Theory 21 (2000), 325-365.
\bibitem{ranick}A. Ranicki, The algebraic construction of the Novikov complex of a circle-valued Morse function, available as math.DG/9903090
\bibitem{schuet}D. Sch\"utz, Gradient flows of closed 1-forms and their closed orbits, available as math.DG/0009055, to appear in Forum Math.
\end{thebibliography}
\end{document}